\documentclass{article}

\usepackage{enumitem}
\usepackage{url}
\usepackage{amsmath}
\usepackage{amssymb}
\usepackage{graphicx}
\usepackage{float}
\usepackage{caption}
\usepackage{subcaption}
\usepackage{amsthm}
\usepackage{color}
\definecolor{green}{rgb}{0,0.5977,0}
\usepackage[linesnumbered,algoruled,boxed,lined]{algorithm2e}
\usepackage[T1]{fontenc}
\usepackage[utf8]{inputenc}
\usepackage{listings}

\definecolor{keywordcolor}{rgb}{0.7, 0.1, 0.1}   
\definecolor{tacticcolor}{rgb}{0.0, 0.1, 0.6}    
\definecolor{commentcolor}{rgb}{0.4, 0.4, 0.4}   
\definecolor{symbolcolor}{rgb}{0.0, 0.1, 0.6}    
\definecolor{sortcolor}{rgb}{0.1, 0.5, 0.1}      
\definecolor{attributecolor}{rgb}{0.7, 0.1, 0.1} 

\newcommand{\zz}{\mathbb Z}

\newcommand{\abs}[1]{\left|{#1}\right|}

\newcommand{\txt}[1]{\text{#1}}

\newcommand{\ignc}[3]{\begin{figure}[H]\begin{center}\includegraphics[scale = {#1}]{#2.png}\caption{#3}\end{center}\end{figure}}

\newcommand{\ipnc}[3]{\begin{figure}[H]\begin{center}\includegraphics[scale = {#1}]{#2.pdf}\caption{#3}\end{center}\end{figure}}

\makeatletter 
\g@addto@macro{\@algocf@init}{\SetKwInOut{Parameter}{Parameters}} 
\makeatother

\newcommand{\lp}[3]{\begin{tabular}{l l}\textbf{#1} & \begin{tabular}{l}$#2$\end{tabular}\\\textbf{subject to} & \begin{tabular}{l l}#3\end{tabular}\end{tabular}}

\newcommand{\lpmin}[2]{\lp{minimize}{#1}{#2}}

\makeatletter
\newcommand*{\shifttext}[2]{%
	\settowidth{\@tempdima}{#2}%
	\makebox[\@tempdima]{\hspace*{#1}#2}%
}
\makeatother

\theoremstyle{plain}
\newtheorem{theorem}{Theorem}
 
\newtheorem{proposition}[theorem]{Proposition}

\newtheorem{conjecture}[theorem]{Conjecture}

\theoremstyle{definition}

\usepackage{tagging}
\usetag{name}

\numberwithin{theorem}{section}

\begin{document}
	\title{Locked Polyomino Tilings}
	\author{\iftagged{name}{Jamie Tucker-Foltz\footnote{\texttt{j.tuckerfoltz@yale.edu}}\\Yale University}{Manuscript \#81}}
	\maketitle
	
	\begin{abstract}
		A locked $t$-omino tiling is a grid tiling by $t$-ominoes such that, if you remove any pair of tiles, the only way to fill in the remaining $2t$ grid cells with $t$-ominoes is to use the same two tiles in the exact same configuration as before. We exclude degenerate cases where there is only one tiling overall due to small dimensions. It is a classic (and straightforward) result that finite grids do not admit locked 2-omino tilings. In this paper, we construct explicit locked $t$-omino tilings for $t \geq 3$ on square and rectangular grids of various dimensions. We show that locked 3- and 4-omino tilings exist on finite square grids of arbitrarily large size. The result for 4-omino tilings in particular is remarkable because they are so rare and difficult to construct: Only a single tiling is known to exist on any grid up to size $40 \times 40$. We also find locked 5- and 6-omino tilings, and provide a computer-assisted proof that locked 7-omino tilings do \emph{not} exist on any rectangular grid. In a weighted version of the problem where vertices of the grid may have weights from the set $\{1, 2\}$ that count toward the total tile size, we demonstrate the existence of locked tilings on arbitrarily large square weighted grids with only 6 tiles. Finally, we construct locked $t$-omino tilings of the infinite grid for arbitrarily large $t$.
		
		Locked $t$-omino tilings arise as obstructions to widely used political redistricting algorithms in a model of redistricting where the underlying census geography is a grid graph. Most prominent is the \emph{ReCom} Markov chain, which takes a random walk on the space of redistricting plans by iteratively merging and splitting pairs of districts (tiles) at a time. Locked $t$-omino tilings are isolated states in the state space of ReCom. The constructions in this paper are counterexamples to the meta-conjecture that ReCom is irreducible on graphs of practical interest.
	\end{abstract}
	\newpage
	\section{Introduction}\label{secIntro}
	
	In the United States, \emph{redistricting} is the process by which states regularly  redraw electoral districts in accordance with new census data. At a technical level, political redistricting is a graph partitioning problem. The vertices represent census blocks, precincts, or counties, and the edges represent geographic adjacency. Given such a vertex-weighted graph $G$ and a target number of districts $k$, redistricting is the act of partitioning the vertex set of $G$ into $k$ sets of equal population, each inducing connected subgraphs.
	
	One of the most successful algorithmic paradigms for exploring the space of redistricting plans (graph partitions) is to run a Markov chain whereby adjacent districts are merged into a double-district and then instantly split again, hopefully in a different way \cite{ReCom}. Such an operation is called a \emph{recombination move}, and the sampling algorithm is known as \emph{ReCom} for short. In practice, the ReCom chain is observed to be rapidly mixing, meaning that, after a short number of steps, it will not matter what redistricting plan the chain started from, and ultimately all plans should have some nonzero probability of being sampled by the algorithm. In theory, however, very little is known about the mixing time.
	
	Arguably the most simple model of this problem is to take $G$ to be an $m \times n$ grid graph, which we notate $G(m, n)$. Grid graphs are simple objects that capture many salient features of those graphs encountered in real-world geographic data. As such, they have been the subject of many prior works on redistricting algorithms \cite{DGKO18, ReCom, RevReCom, procacciaTuckerFoltz2021compactness, tappspanning, charikar2022complexity, frieze2022subexponential, TreeSplitting, NnToN}. We assume, as is typical in this model, that vertices have equal population. (This assumption is revisited in Section~\ref{secSixDistricts}.) The \emph{metagraph} $\mathcal{M}(G, t)$ is the graph whose vertices are partitions of $G$ into connected subgraphs, each containing exactly $t$ vertices, with two partitions joined by an edge whenever they differ only on two parts (see Figure \ref{figMetagraphExample}). A ReCom Markov chain on redistrictings of $G$ into $k$ districts is thus a walk on $\mathcal{M}(G, t)$, where $t = \frac{\abs{G}}{k}$, so, given a probabilistic transition function, we may inquire about its mixing time. A prerequisite, of course, is the basic question of whether or not $\mathcal{M}(G, t)$ is even connected. Unfortunately, even for $k = 3$, very little is known.
	
	\ipnc{.16}{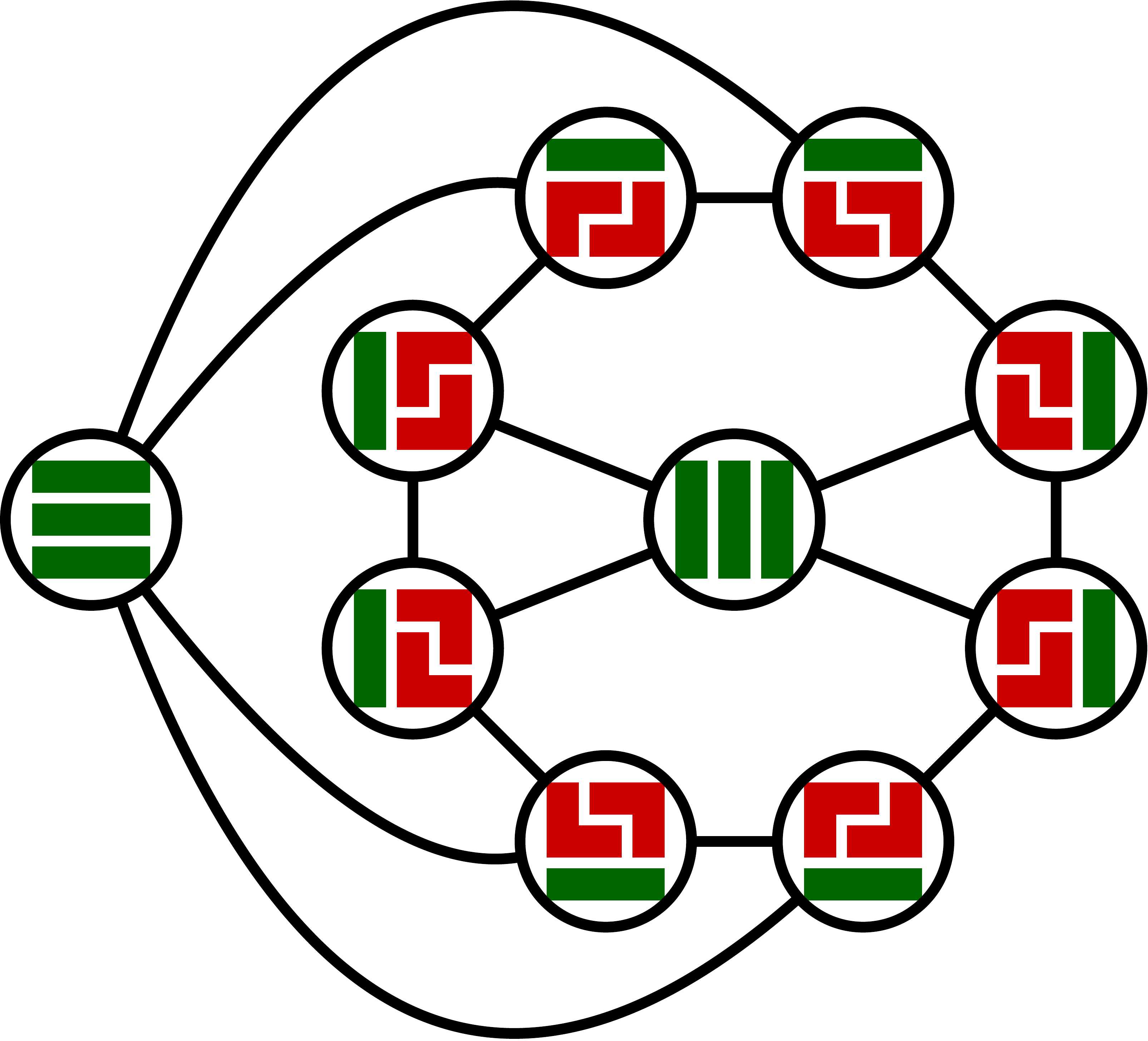}{\label{figMetagraphExample}The metagraph $\mathcal{M}(G(3, 3), 3)$. Each vertex is a partition of the $3 \times 3$ grid graph into connected subgraphs, each containing 3 vertices. Note that this graph is connected, as is claimed by Conjecture \ref{cnjGridKEquals3}.}
	
	\begin{conjecture}\label{cnjGridKEquals3}
		For any positive integer $n$ that it divisible by 3, $\mathcal{M}(G(n, n), \frac{n^2}{3})$ is connected.
	\end{conjecture}
	
	Every step, you are allowed to repartition two-thirds of the entire grid graph any way you wish, so surely it is possible to eventually get from one partition to any other? A generalization\footnote{Note that Conjecture~\ref{cnjGridArbitraryK} is not a strict generalization of Conjecture~\ref{cnjGridKEquals3}, since holding ``for all sufficiently large $n$'' is a weaker condition. However, Conjecture~\ref{cnjGridKEquals3} holds for all small grid sizes for which the metagraph has been fully enumerated, and it would be quite strange if there was a large counterexample that does not generalize to arbitrarily large grids. Indeed, the counterexamples presented in this paper are all in the regime where $n$ is much smaller than $k$, which is nonexistent for $k = 3$.} of this conjecture to any $k$ is as follows.
	
	\begin{conjecture}\label{cnjGridArbitraryK}
		For any positive integer $k$, for all sufficiently large $n$ such that $k$ divides $n^2$, $\mathcal{M}(G(n, n), \frac{n^2}{k})$ is connected.
	\end{conjecture}
	
	If we completely drop the constraint that parts have equal numbers of vertices, it is known that the metagraph of partitions into $k$ districts is connected whenever $G$ is 2-connected, and moreover, it is possible to transition from one partition to another by only reassigning a single vertex at each step \cite{AnySizeDistricts}. When $G$ has a Hamiltonian cycle, Akitaya, Korman, Korten, Souvaine, and T\'oth \cite{AnySizeDistrictsNew} show that districts need be no bigger than double their target size (but may still need to be as small as a single vertex). More relevant is the groundbreaking triangular lattice metagraph connectivity result for $k = 3$ by Cannon \cite{TriangularLattice}, which keeps very strict balance constraints. Compared to Conjecture \ref{cnjGridKEquals3}, this result is weaker in the following three ways.
	\begin{itemize}
		\item The graph $G$ is a subset of the triangular lattice rather than a square lattice. This graph structure is easier to work with since the neighbors of any vertex are connected by a path/cycle, a fact that is leveraged by the proof in a crucial way.
		\item The districts may vary in size $\pm 1$ from the target size. This makes it easier to transition from one partition to another, as only one vertex need be reassigned in a given step.
		\item Districts must always be simply-connected.
	\end{itemize}
	
	Subsequent work by Akitaya, Dituro, Gonczi, Korman, Souvaine, Stock, and T\'oth~\cite{SquareGridsSosa} managed to simplify this proof considerably, and also extend to the square lattice under the same conditions.
	
	In this paper, we approach the metagraph connectivity question from the opposite end of the parameter space: instead of taking $k$ to be small, we assume the district size $t$ to be small. On grid graphs, it is natural to think of this problem in the language of \emph{polyomino tilings}. Specifically, we study the space of $t$-omino tilings of $m \times n$ grids under operations that replace two tiles at a time. Here we find new obstacles, exemplified by the 3-omino tiling in Figure \ref{fig3Omino6}.
	
	\ipnc{.16}{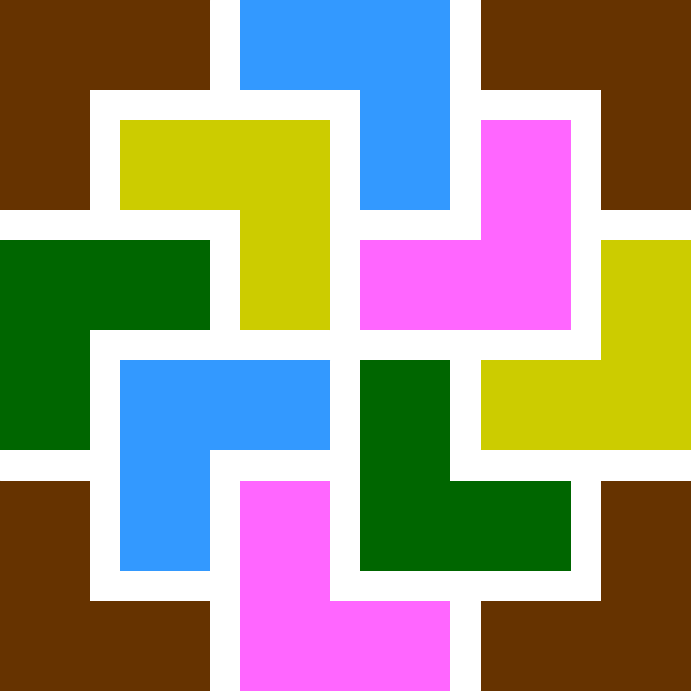}{\label{fig3Omino6}A locked 3-omino tiling of $G(6, 6)$. In total, there are 50 such tilings.}
	
	Observe that there are no two adjacent tiles that can be merged and split into 3-ominoes in a different way than they are currently arranged. Since there are no recombination moves, this tiling is an isolated vertex in $\mathcal{M}(G(6, 6), 3)$, proving that the graph is disconnected. Further, observe that we may tile any $6m \times 6n$ grid with copies of this tiling arranged in an $m \times n$ grid of $6 \times 6$ blocks as in Figure \ref{fig3Omino12}. By flipping every other block, we ensure that there are no recombination moves available on the boundaries between blocks either. This proves that grid metagraphs can be disconnected for arbitrarily large $k$, even fixing $t = 3$. In view of Conjecture~\ref{cnjGridArbitraryK}, this shows that requiring the statement to only hold for ``sufficiently large $n$'' is necessary.
	
	Formally, we define a \emph{locked $t$-omino tiling} of an $n$-vertex graph $G$ to be an isolated vertex in $\mathcal{M}(G, t)$ that is not the only vertex in $\mathcal{M}(G, t)$. This final condition outlaws degenerate cases: a tiling cannot be considered locked merely because there are no other tilings (as this is not a connectivity counterexample anyway). The fundamental question we ask in this paper is,
	
	\ipnc{.16}{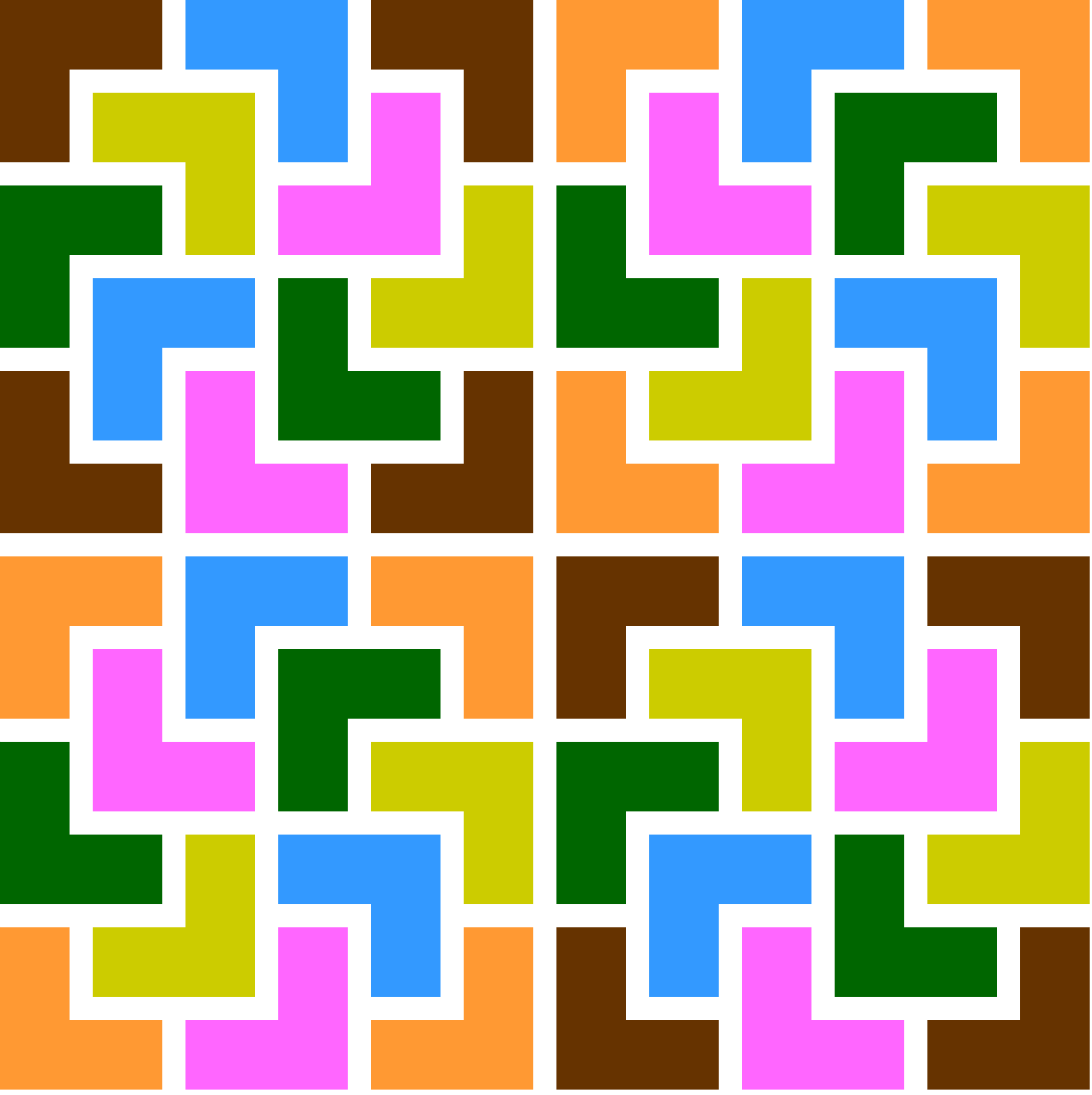}{\label{fig3Omino12}A locked 3-omino tiling of $G(12, 12)$ obtained from copying and flipping the tiling from Figure \ref{fig3Omino6}.}
	
	\begin{center}
		\emph{For which grid dimensions and which $t$ do there exist locked $t$-omino tilings?}
	\end{center}
	This question is interesting for two principal reasons:
	\begin{enumerate}
		\item Locked tilings point to the theoretical limits of ReCom and other MCMC-based redistricting algorithms. Much effort has been spent on trying to prove rapid mixing. Understanding the combinatorial obstructions that can occur informs what assumptions we might need to make in order to guarantee mixing properties, or what tweaks to the algorithms may be necessary (e.g., relaxing balance constraints, recombining more than two districts at a time). Furthermore, the problem of balanced graph partitioning with small parts arises when jurisdictions impose \emph{nesting rules}, for example, that a state senate district must be formed from exactly 3 state house districts \cite{NestingRules}.
		\item There is a substantial mathematics literature on connectivity of polyomino tilings via local moves. This particular question, where we do not restrict the set of $t$-ominoes that may be used, is arguably the most basic one to ask, yet it has not been studied before beyond the case of $t = 2$ (see Section \ref{subLitReview}). Perhaps this is due to the fact that locked $t$-omino tilings for larger $t$ are so difficult to find that they require computational methods, as this paper employs (see Section \ref{subAlgorithm}).
	\end{enumerate}
	
	\begin{table}[H]\begin{center}
			\begin{tabular}{ r|l }
				Polyomino size $t$ & Existence on rectangular grids\\
				\hline
				2 & Do not exist\\
				3 & Exist on all rectangular grids of sufficiently large dimensions\\
				4 & Exist on infinitely many square grids\\
				5 & Exist, and thousands have been found\\
				6 & Exist, but only one has been found\\
				7 & Do not exist
			\end{tabular}
			\caption[]{\label{tabExistenceSummary}For which $t$ are there any locked $t$-omino tilings of rectangular grids.}
	\end{center}\end{table}
	
	Table \ref{tabExistenceSummary} summarizes the landscape of our results for finite rectangular grids and constant $t$. Prior to this work, only the first row was known: For $t = 2$, no locked polyomino tilings exist on any rectangular grid. In the special case of the $8 \times 8$ grid, this is a classic mathematical puzzle: Is it possible to cover a chessboard with $1 \times 2$ dominoes without creating a $2 \times 2$ square? Not only is this not possible, but it is well-known that the metagraph is connected on any rectangular grid, and moreover, ReCom is rapidly mixing \cite{GlauberConnected2, GlauberConnected1}. For completeness, in Section \ref{sec2Omino} we give an elementary proof of the connectivity of 2-omino metagraphs on grids.
	
	In Section \ref{sec3Omino} we then turn to 3-ominoes, where locked tilings are abundant. The example from Figures \ref{fig3Omino6} and \ref{fig3Omino12} show that grids with side lengths divisible by 6 admit locked tilings; we prove a stronger result, that all sufficiently large grids admit locked tilings. On the other hand, if one of the side lengths is at most 3, we show that the metagraph is connected, so no locked tilings exist.
	
	In Section \ref{sec45Omino} we construct locked 4-, 5-, and 6-omino tilings, which are considerably more difficult to find. The constructions in this section were discovered in part with an extensive computational search. One would expect to find none, as is the case with 2-ominoes, or many, as is the case with 3-ominoes. Instead, after exhaustively searching for all locked tilings on the $m \times n$ grid for any $m, n \leq 20$, the program found only a \emph{single} 4-omino and 5-omino tiling up to isomorphism: There is one locked 4-omino tiling on a $10 \times 10$ grid, and one locked 5-omino tiling on a $20 \times 20$ grid. The next smallest known locked 4- and 5-omino tilings are on the $42 \times 42$ and $30 \times 30$ grids, respectively, while the only known 6-omino tiling is on a $48 \times 48$ grid. In Section \ref{sec45Omino} we describe the algorithm and present these tilings, along with a generalization of the $42 \times 42$ tiling to an infinite family of locked 4-omino tilings on larger and larger finite grids. We also present a computer-assisted proof that there do not exist any locked 7-omino tilings, which has been formalized in Lean.
	
	In Section \ref{secTori} we briefly consider infinite periodic tilings, in both the square and triangular grid lattices. These are easier to reason about since there is no boundary, so all vertices look the same. While we conjecture that the set of integers $t$ for which there exist locked $t$-omino tilings of finite grids is bounded (possibly even by 5), it turns out that this is not the case for infinite tilings. We construct a family of infinite periodic locked $t$-omino tilings for increasing $t$. We also extend this result to the triangular lattice.
	
	In Section \ref{secSixDistricts} we consider a seemingly innocuous alteration of our grid graph model: Vertices may have population weight 1 or 2. Clearly, if we allowed arbitrary integral population weights, metagraph connectivity is hopeless, as getting exact balance could require satisfying arbitrarily complex constraints about certain vertices being in the same parts. It is reasonable to expect that we cannot construct such pathological examples with bounded-size population fluctuations. However, in Theorem~\ref{thmSixDistricts}, we do just that. We demonstrate a family of arbitrarily large instances on square grid graphs on which the metagraph is not connected for $k = 6$. In contrast to all of the other results of our paper, this result holds in a regime with a constant number of districts of a large size. It shows that Conjecture~\ref{cnjGridArbitraryK} cannot be generalized beyond grid graphs to even extremely tame families of ``grid-like graphs'' that have been studied in prior work on redistricting algorithms.
	
	\subsection{Related Work on Polyomino Recombination}\label{subLitReview}
	
	The study of connectivity of domino tilings by local moves goes back to the works of Conway and Lagarias \cite{Conway} and Thurston \cite{GlauberConnected1}, who construct a ``height function'' on the dual graph with the property that the highest point can always be decreased by a recombination move. This actually turns the set of domino tilings into a lattice, implying that the state space is connected \cite{ThurstonMadePrecise}. Surprisingly, this holds not just on rectangular grids, but any any simply connected grid subgraph. Reachability among infinite domino tilings has also been studied \cite{WholePlane}; unsurprisingly, not all pairs are connected, even by an infinite number of moves.
	
	Domino tilings generalize to higher dimensions as $2 \times 1 \times 1 \times \dots \times 1$ blocks. Similar questions have been asked in finite grids and grid subgraphs of dimension three \cite{ThreeDimensionalDominoes1, ThreeDimensionalDominoes2} and larger \cite{FourDimensionalDominoes}. Even on the rectangular $3 \times 3 \times 2$ grid, domino tilings many be disconnected under pairwise recombination, leading to the study of three-way ``trit'' moves.
	
	The space of domino tilings can also fail to be connected on two-dimensional lattices other than the square lattice. If we add diagonals to the square lattice (all in the same direction), we get the triangular lattice, on which there are simple examples of locked tilings. For instance, consider a triangular region with six vertices. There are only two perfect matchings in the triangular lattice, and they are not connected by a recombination move \cite{TriangularLattice}. However, Kenyon and R\'emila \cite{DominoesTriangularLattice} showed that four-way recombination moves suffice to connect the space of domino tilings on any simply connected subgraph of the triangular lattice. R{\o}ising and Zhang \cite{DominoesArchimedianLattices} generalized this result to all eleven Archimedean lattices. If we add diagonals to the square grid in \emph{alternating} directions, then pairwise recombination moves suffice \cite{Impurities}.
	
	Larger polyomino sizes have received far less attention. To the best of our knowledge, this is the first paper to study the question of state space connectivity with arbitrary tiles of a given constant size $t \geq 3$. Related questions have been investigated previously in the context of restricted tile sets. Korn \cite{KornThesis} studies the question of when arbitrarily powerful local moves, perhaps involving recombining more than two tiles, exist for restricted sets of polyominoes. While not every tile set admits \emph{any} set local moves, it is known that recombinations of up to 4 pieces are sufficient to connect the space of grid tilings by T-shaped 4-ominoes \cite{TOmino}. Sheffield \cite{Ribbon} shows that pairwise recombination moves connect the set of \emph{ribbon tilings} of order $n$ in any simply connected grid subgraph. In a ribbon tiling, each piece must be a path that snakes upward and to the right. Pak~\cite{PakSurvey} surveys several other results about recombination moves and invariants for finite restricted tile sets.
	
	\section{Domino Tilings}\label{sec2Omino}
	
	In this section we present an elementary proof of the following result, which is originally credited to William Thurston \cite{GlauberConnected1} as mentioned above. Our proof, while less elegant and general, is perhaps easier to understand.
	
	\begin{theorem}\label{thmDominoConnected}
		For any positive integers $m$ and $n$, $\mathcal{M}(G(m, n), 2)$ is connected.
	\end{theorem}

	\begin{proof}
		If either $m$ or $n$ is 1, or if $mn$ is odd, then there is at most one tiling, so the metagraph is obviously connected. So suppose both side lengths are at least 2, and without loss of generality assume the horizontal dimension has even length. We show that any pair of tilings are connected to each other via a sequence of recombination moves by showing that any given tiling is connected to a ``central'' tiling in which all dominoes are horizontally oriented. Starting from the top-left corner, reading left-to-right, top-to-bottom, let $D_1$ be the first domino that is vertically oriented. In Figure \ref{figDominoConnectivityProof}, $D_1$ is the black-outlined domino.
		
		\ipnc{.16}{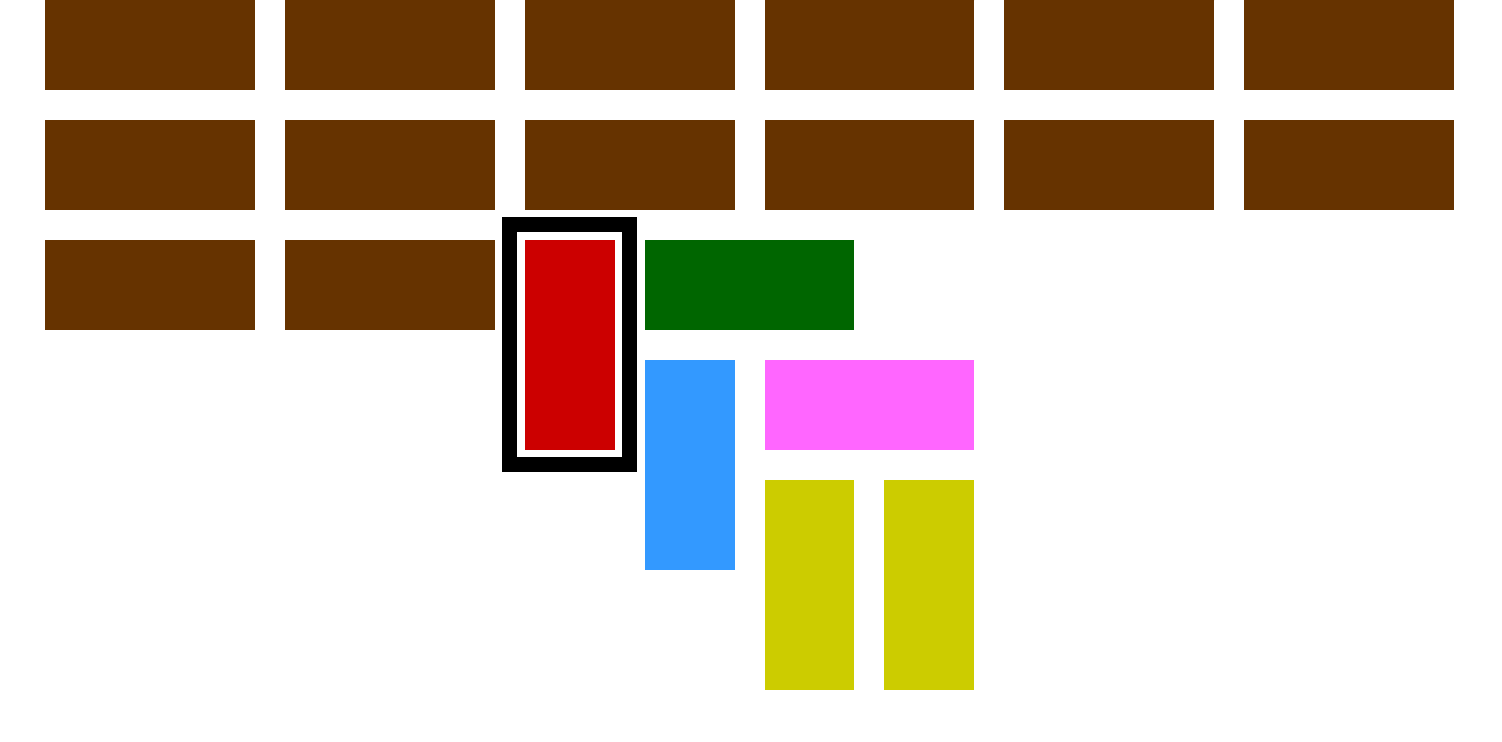}{\label{figDominoConnectivityProof}A generic step of the algorithm to transition to the all-horizontal domino tiling.}
		
		If the domino $D_2$ to the right of $D_1$ is also vertically oriented, then we can recombine $D_1$ and $D_2$ so that $D_1$ is horizontal. Otherwise we check if we can recombine $D_2$ with the domino $D_3$ below it, and so on. Iteratively supposing we cannot recombine, we find a ladder of dominoes extending downward and to the right, which must eventually terminate in a $2 \times 2$ square (the gold vertical dominoes in the bottom of Figure \ref{figDominoConnectivityProof}) because it cannot go on forever without hitting the boundary of the rectangular grid. It is straightforward to see that a sequence of recombination moves can push this square up the ladder until $D_1$ is horizontal. Note that these moves do not touch any of the horizontal dominoes preceding $D_1$. Thus, we have strictly increased the number of consecutive horizontal dominoes reading from the top-left corner. We can continue iteratively applying this algorithm until all dominoes are horizontal.
	\end{proof}
	
	\section{Tromino Tilings}\label{sec3Omino}
	
	We begin our discussion of 3-ominoes by considering the easy cases where the metagraph is connected. The proof proceeds along similar lines as the proof of Theorem \ref{thmDominoConnected}, though there are a few more cases to consider.
	
	\begin{theorem}\label{thmTriominoConnected}
		For any positive integers $m$ and $n$ such that $m \leq 3$, $\mathcal{M}(G(m, n), 3)$ is connected.
	\end{theorem}
	
	\begin{proof}
		As before, we may assume without loss of generality that $mn$ is divisible by 3 and both $m$ and $n$ are at least 2.
		
		\ipnc{.16}{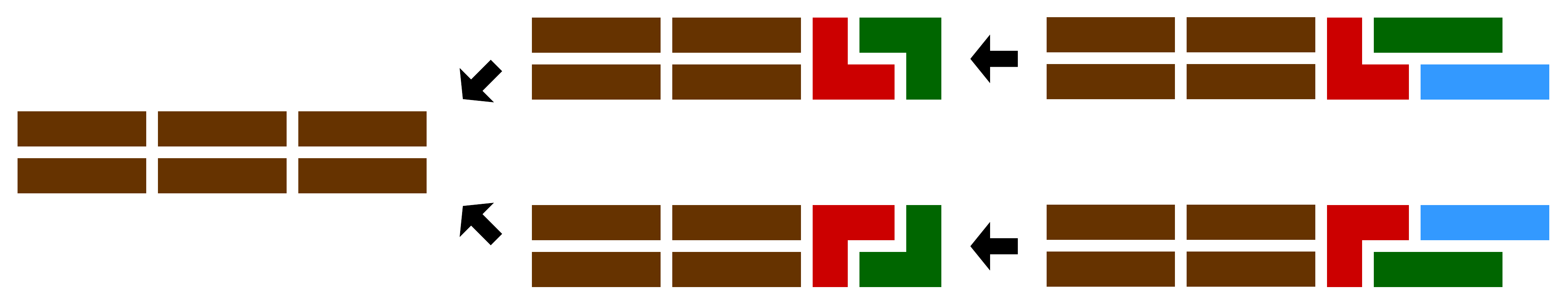}{\label{figTriominoConnectivityProof2} Steps required to make all 3-ominoes horizontal in a $2 \times n$ grid.}
		
		If $m = 2$, then we show that we can transition from any tiling to a tiling in which all 3-ominoes are horizontal. Indeed, the first column that is not configured this way must be configured according to one of the four cases in Figure \ref{figTriominoConnectivityProof2}. As indicated by the arrows, we can reconfigure these 3-ominoes to be horizontal in at most 2 moves. Continuing inductively, we can eventually make all 3-ominoes horizontal.
		
		\ipnc{.16}{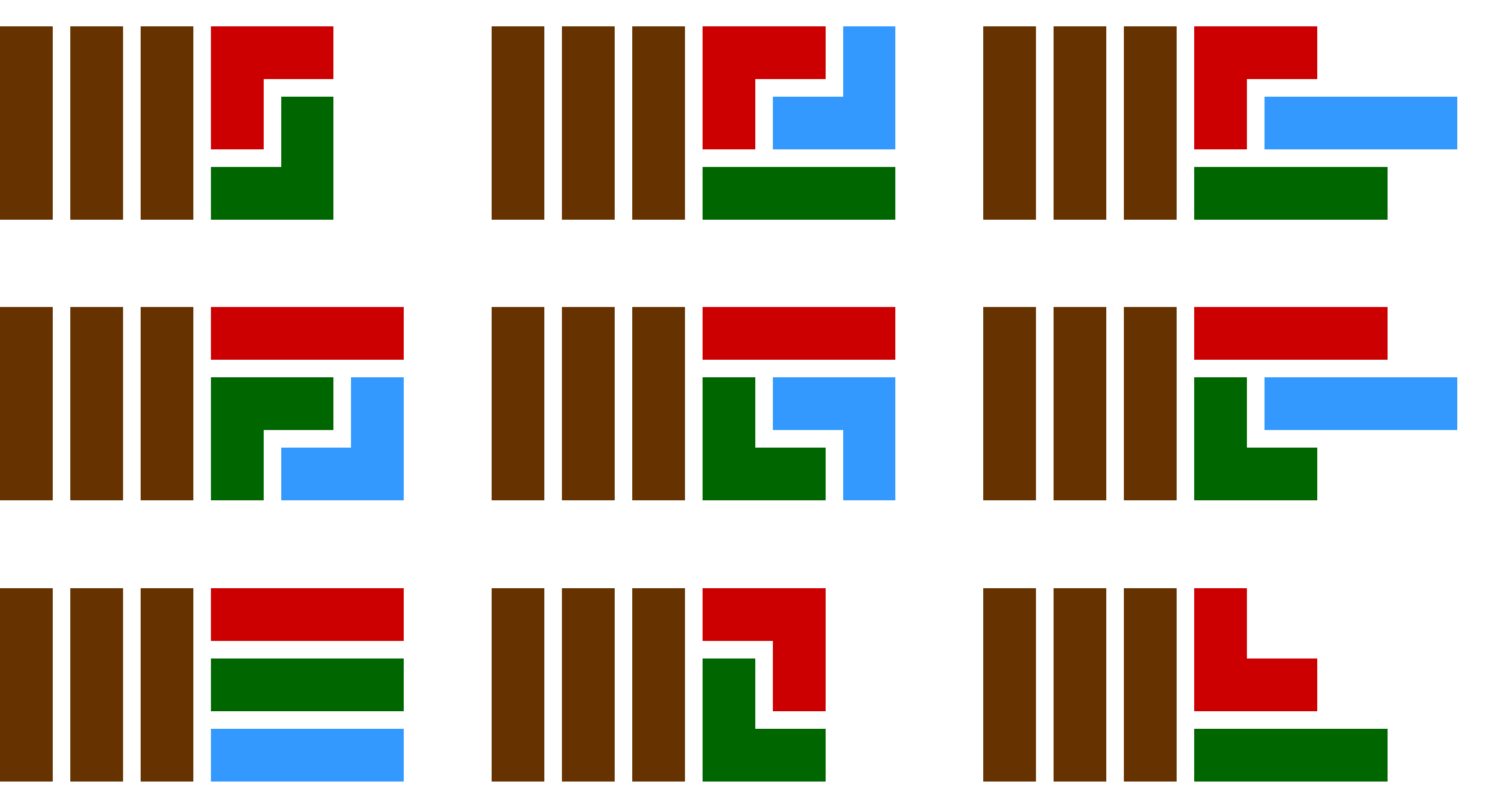}{\label{figTriominoConnectivityProof3} All possible cases of what the frontier could look like in transitioning to an all-vertical tiling of a $3 \times n$ grid.}
		
		If $m = 3$, we analogously transition to a tiling in which all 3-ominoes are \emph{vertical}. Figure \ref{figTriominoConnectivityProof3} exhaustively lists all possible cases of what the neighborhood around the first column that is not a 3-omino could look like. We leave it to the reader to verify that it is always possible to make this column into a vertical 3-omino in at most 3 recombination moves.
	\end{proof}
	
	The smallest nontrivial case that is not covered by this theorem is $m = 6$ and $n = 4$. Already here there are locked 3-omino tilings---up to isomorphism, there are 7 distinct tilings. We conjecture that locked tilings exist whenever Theorem \ref{thmTriominoConnected} does not apply, i.e., whenever $m$ is divisible by 3 (which is without loss of generality), $m \geq 6$, and $n \geq 4$. It is tedious to prove a statement like this as, even though 3-omino tilings are plentiful, they are largely unstructured and difficult to find, especially when $m$ or $n$ is odd. With not to much difficulty, we can obtain a weaker result, replacing the number 4 with 100. This at least shows that there are no grid dimension parity obstructions to producing locked tilings.
	
	\begin{theorem}\label{thm3OminoLocked}
		For any integers $m$ and $n$ such that $m$ is divisible by 3, $m \geq 6$, and $n \geq 100$, the $m \times n$ grid has a locked 3-omino tiling.
	\end{theorem}
	
	\begin{proof}
		Given the assumptions on $m$ and $n$, one can verify that it is always possible to find nonnegative integers $x_1, x_2, x_3, x_4, x_5$ such that $m = 6x_1 + 9x_2$ and $n = 22x_3 + 19x_4 + 8x_5$, with one of $x_1$ or $x_2$ greater than zero and $x_3$ also greater than zero.\footnote{That such values exist for large enough $m$ and $n$ follows from the facts that $\gcd(6, 9) = 3$ and $\gcd(22, 19, 8) = 1$. It is a straightforward and tedious calculation to check that $m \geq 6$ and $n \geq 100$ are large enough thresholds.} We may thus draw straight horizontal and vertical lines across the entire grid to partition the grid into rectangular blocks of size $6 \times 22$, $9 \times 22$, $6 \times (22 + 19x_4 + 8x_5)$, and $9 \times (22 + 19x_4 + 8x_5)$. We may fill these blocks in with the tilings depicted in Figure \ref{figTriominoLocked}.
		
		\ipnc{.16}{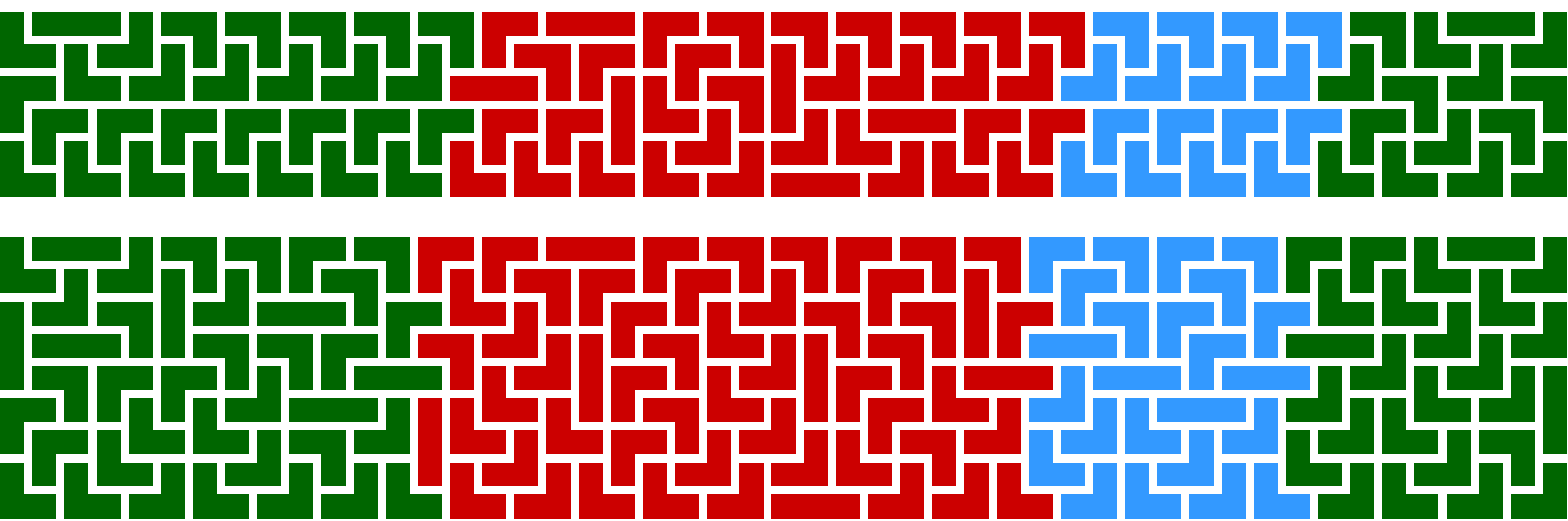}{\label{figTriominoLocked} Locked 3-omino tilings that can tile with themselves and each other.}
		
		Each tiling consists of four parts: two green parts on the end of total length 22, a red part in the middle of length 19, and a blue part of length 8. To make a rectangle of length $22 + 19x_4 + 8x_5$ we repeat the red part $x_4$ times and the blue part $x_5$ times. One can check that there are no recombination moves within the interiors of any of these tilings, nor between neighboring blocks.
	\end{proof}
	
	\section{Computational Findings for $t \geq 4$}\label{sec45Omino}
	
	\subsection{Enumeration Algorithm}\label{subAlgorithm}
	
	While locked 3-omino tilings can be found by hand with some effort, locked $t$-omino tilings are extremely difficult to construct for $t \geq 4$. However, for small grids we may attempt to enumerate them exhaustively using a computer. The surprising takeaway from these computational results is that locked 4- and 5-omino tilings exist, but are extraordinarily rare.
	
	The algorithm, which has been made publicly available online,\footnote{\iftagged{name}{\url{https://github.com/jtuckerfoltz/LockedPolyominoTilings}}{Link excluded for the purposes of anonymity.}} works as follows. We enumerate all possible locations within all possible $t$-ominoes that a given cell could be contained in, which we call the \emph{type} of the cell. For example, when $t = 2$, there are 4 types, depending on whether the other cell in the domino is above, below, to the left, or to the right. For $t = 3$, there are 18 types, for $t = 4$, there are 76 types, and for $t = 5$, there are 315 types.\footnote{See OEIS sequence A048664 \cite{OtherOEIS}.} Next, we enumerate, for every $(x, y)$ offset, all the incompatible pairs of types, i.e., the set of pairs of types $(t_1, t_2)$ such that it is not possible for a cell at location $(x_0, y_0)$ to have type $t_1$ and a cell at location $(x_0 + x, y_0 + y)$ to have type $t_2$ because the tiles either intersect or admit a recombination move. We store the incompatible pairs in a massive lookup table, representing types as arbitrary integers to save space. We then augment the lookup table by iteratively applying the following rule until no more additions are made: it is not possible for cell $a$ to have type $t_1$ and cell $b$ to have type $t_2$ if that would mean that there is no compatible type for some other cell $c$. All of these steps are very time consuming, but only need to be computed once for each polyomino size.
	
	For any given grid size, we initialize each grid cell with the set of types it is allowed to have. For most interior cells, this is almost all types, but cells near the grid boundary have additional restrictions. We then iteratively pick a cell with few possible types and recursively try all of them, eliminating incompatible types for other cells as we go. This enables us to exhaustively search for all locked polyomino tilings. The reason this algorithm is so successful at handling large grid sizes is that the locking constraint is quite strong, which severely restricts the branching factor.
	
	\subsection{Computational Results on Small Grids}\label{subSmallGrids}
	
	Using the algorithm described in the previous section, we may enumerate all locked 4-, 5-, and 6-omino tilings on rectangular grids up to size $20 \times 20$, and all tilings that are also rotationally symmetric on grids up to size $40 \times 40$. Table \ref{tab456} summarizes the results: Only a small handful of locked 4- and 5-omino tilings exist, and no locked 6-omino tilings exist.
	
	\begin{table}[ht]\begin{center}
			\begin{tabular}{ r|c|c|c }
				Number of locked $t$-omino tilings on grids of size\dots & 4-ominoes & 5-ominoes & 6-ominoes\\
				\hline
				\dots up to $20 \times 20$ & 1 & 1 & 0\\
				\dots up to $40 \times 40$, with 4-fold rotational symmetry & 1 & 3270 & 0
			\end{tabular}
			\caption[]{\label{tab456}Enumeration of locked 4-, 5-, and 6-ominoes on small grid sizes, counting up to isomorphism.}
	\end{center}\end{table}
	
	
	
\begin{figure}[H]
		\centering
		\begin{minipage}{.5\textwidth}
			\centering
			\includegraphics[scale=.16]{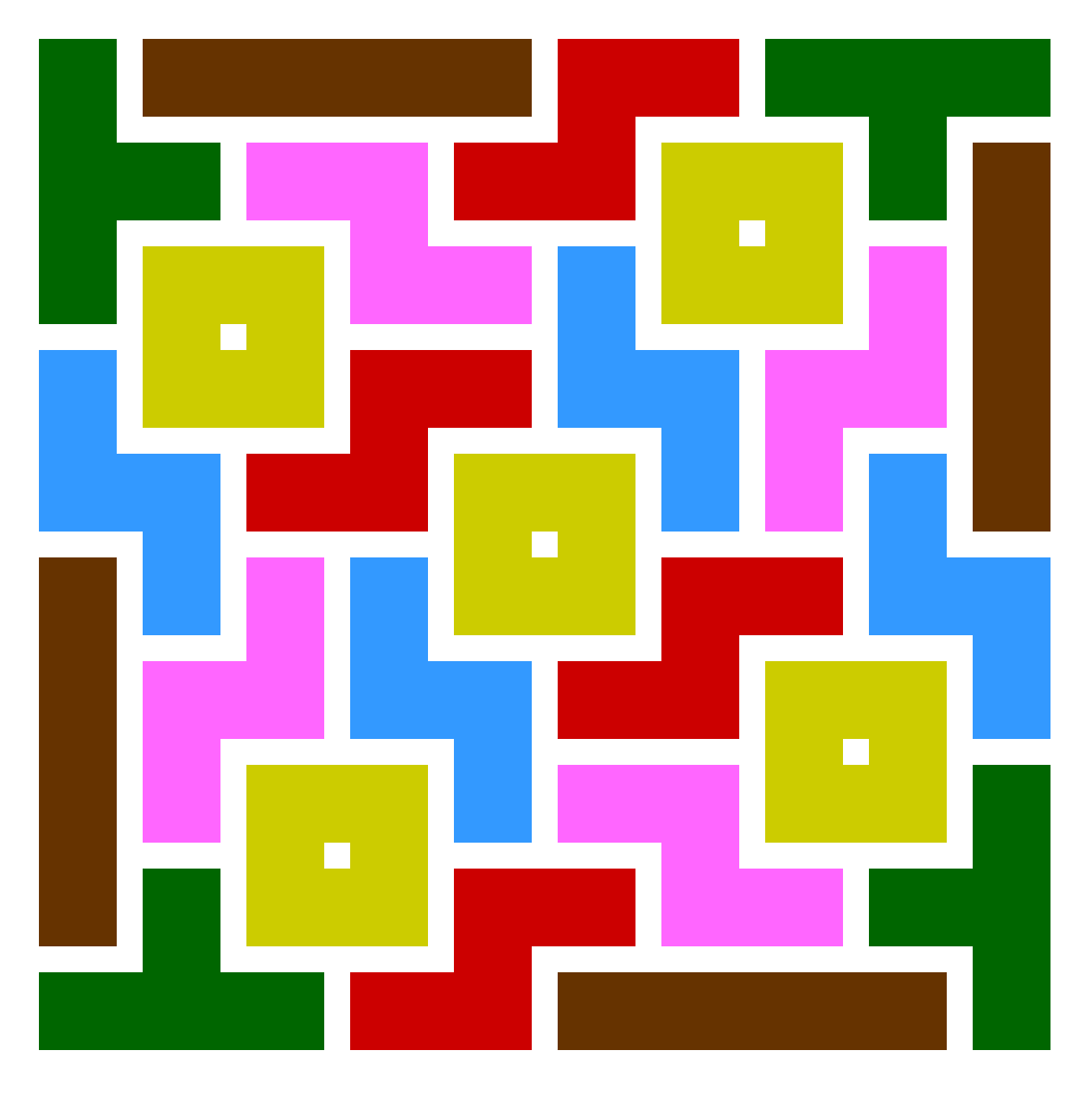}
			\caption{\label{figUniqueTilings10}The smallest locked 4-omino grid tiling.}
		\end{minipage}%
		\begin{minipage}{.5\textwidth}
			\centering
			\includegraphics[scale=.16]{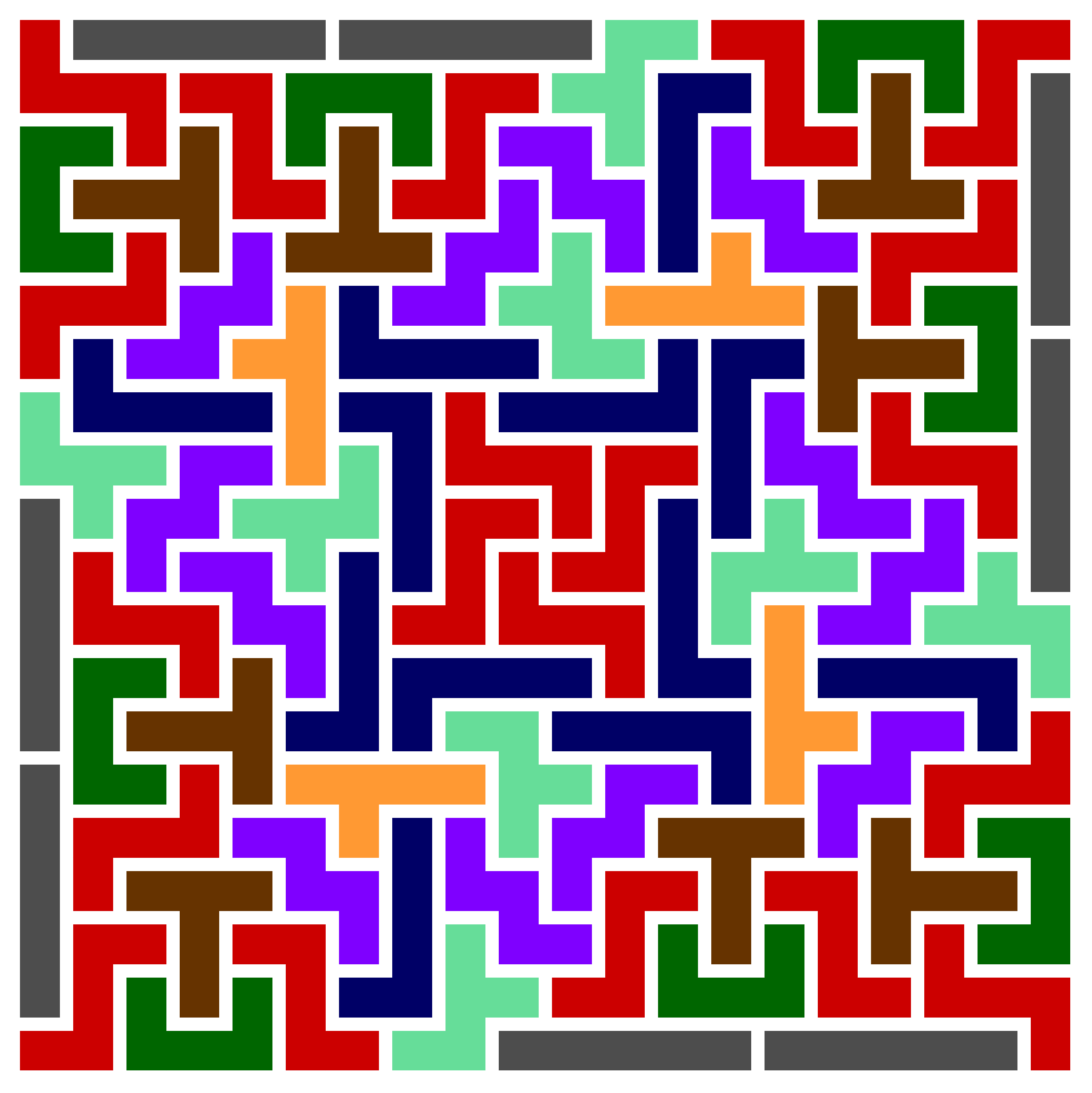}
			\caption{\label{figUniqueTilings20}The smallest locked 5-omino grid tiling.}
		\end{minipage}
	\end{figure}

	
	Figure \ref{figUniqueTilings10} shows the sole locked 4-omino tiling mentioned in this table, and Figure \ref{figUniqueTilings20} shows the smallest locked 5-omino tiling.\tagged{name}{\footnote{Pictures of additional tilings can be found online at:\\\url{https://github.com/jtuckerfoltz/LockedPolyominoTilings/tree/main/Images}}} Table \ref{tabSymmetric5OminoTilings} further breaks down the enumeration of locked 5-omino tilings by grid dimensions.
	
	\begin{table}[H]\begin{center}
			\begin{tabular}{ c|c }
				Square grid side length & Number of symmetric locked 5-omino tilings (A365862)\\
				\hline
				5 & 0\\
				10 & 0\\
				15 & 0\\
				20 & 1\\
				25 & 0\\
				30 & 10\\
				35 & 27\\
				40 & 3232
			\end{tabular}
			\caption[]{\label{tabSymmetric5OminoTilings}Number of 4-fold symmetric locked 5-omino tilings, up to isomorphism. This sequence appears doubled in OEIS A365862 \cite{MyOEIS}, since each tiling can be mirrored to produce an isomorphic one. The sole tiling of the $20 \times 20$ grid is depicted in Figure \ref{figUniqueTilings20}.}
	\end{center}\end{table}

		%
		%

	We record several additional facts about the structure of locked tilings that were discovered computationally.
	
	\begin{enumerate}
		\item Any locked 4-omino tiling of a grid where both dimensions are at least 14 must contain all five tile shapes. (It is notable that the $10 \times 10$ tiling from Figure \ref{figUniqueTilings10} does not; it is missing the L tile.)
		\item While we cannot enumerate the asymmetric locked 5-omino tilings of some larger grid sizes due to computational issues, we do know that at least some tilings exist. In particular, Figure \ref{figSpaceInvaders} shows a partial tiling with four holes that can be completed in different ways to produce asymmetric tilings. They can also be completed to yield small metagraph components of size greater than one.
		\item Any tile can be part of a locked 5-omino grid tiling except for the V-tile. It is not too hard to see that any tile placed inside the V-tile will admit a recombination move with it. There are locked tilings that use all of the other tiles (for instance, those in Figure \ref{figSpaceInvaders}).
		\item In any locked 4-, 5-, or 6-omino tiling of a grid where both dimensions are at least 9, only some specific tiles can be placed in the corners. For the top left corner, the tile must be one of the tiles in Figure \ref{figCorners} (or flipped along the main diagonal).
	\end{enumerate}
	
	\ipnc{.16}{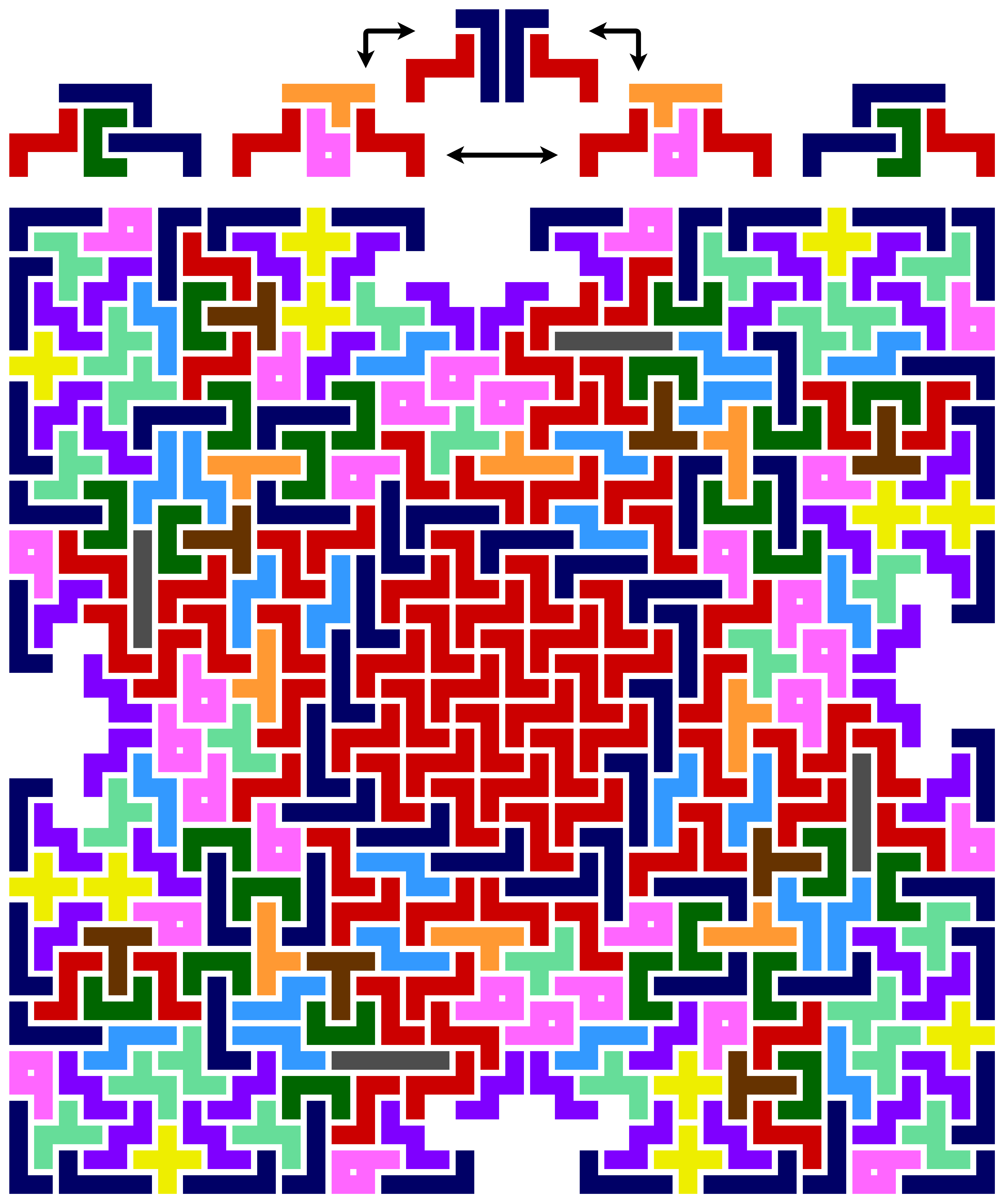}{\label{figSpaceInvaders}A partial tiling of the $40 \times 40$ grid graph. The remaining spaces can each be completed in five different configurations. Two of these configurations are locked and the other three are linked to each other by recombination moves, but can never recombine with any other neighboring tiles. Thus, the tiling can be completed to obtain both symmetric and asymmetric locked tilings, as well as connected components of $\mathcal{M}(G(40, 40), 5)$ of size 3, 9, 27, and 81.}
	
	\ipnc{.16}{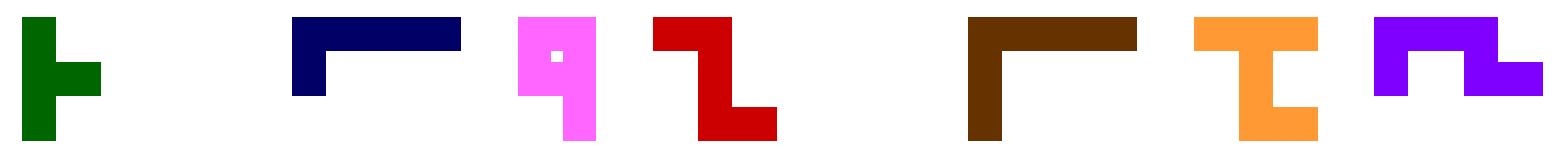}{\label{figCorners}The only $t$-ominoes that can occur in the top-left corner of a locked $t$-omino grid tiling, for $t \in \{4, 5, 6\}$.}
	
	\ipnc{.16}{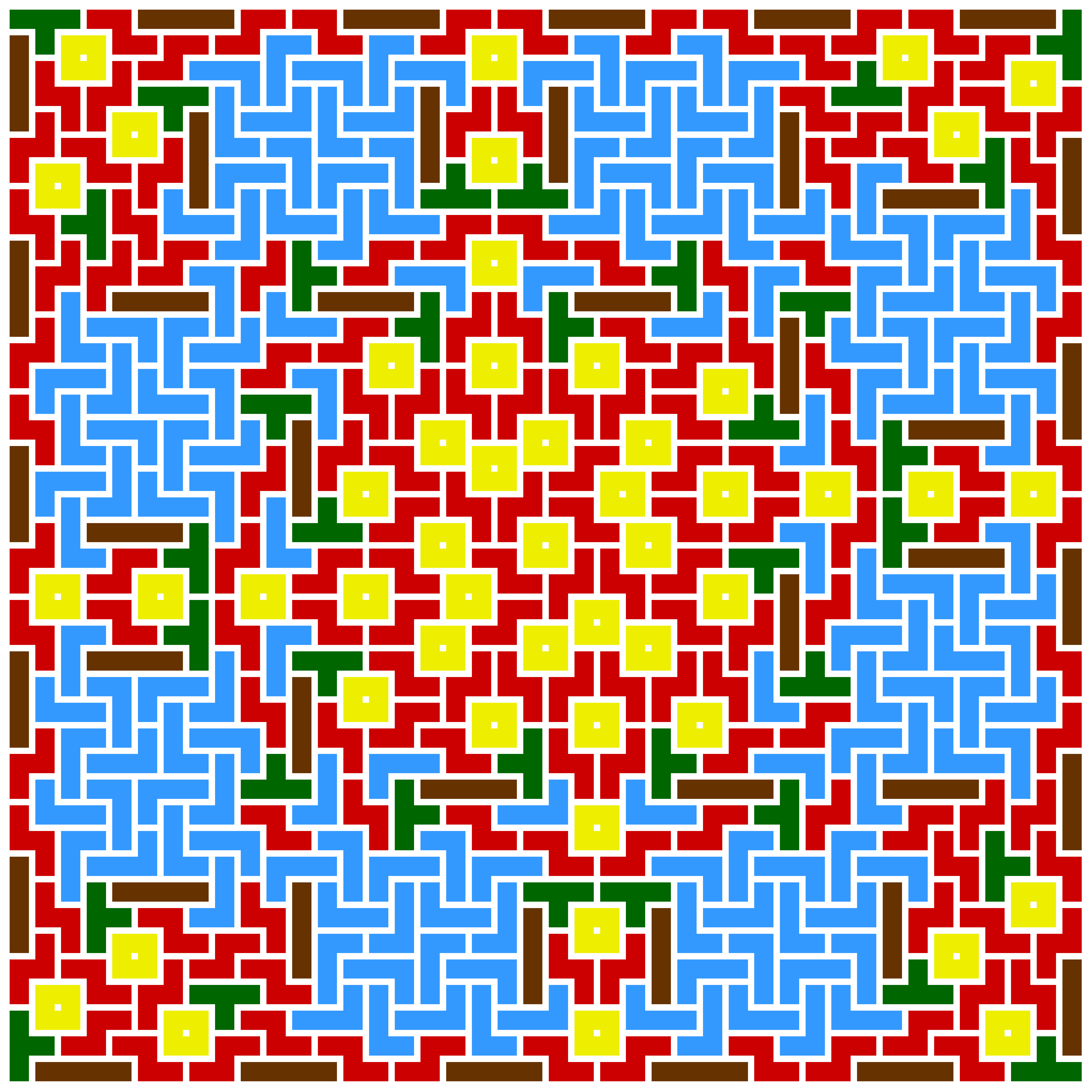}{\label{fig4OminoLockedInfiniteConstructionSmall}A locked 4-omino tiling on the $42 \times 42$ grid. This is the smallest known locked 4-omino tiling beyond the $10 \times 10$ grid tiling from Figure \ref{figUniqueTilings10}. It is the tiling for $n = 1$ in construction of Theorem \ref{thm4OminoLockedInfiniteConstruction}.}
	
	\subsection{Locked 4-Omino Tilings on Arbitrarily Large Grids}
	
	One might conjecture, on the basis of Table \ref{tab456}, that the locked 4-omino tiling from Figure \ref{figUniqueTilings10} is the only locked tiling that exists on any finite grid. Surprisingly, not only do more tilings exist, but they exist for arbitrarily large grid dimensions, starting at $42 \times 42$.
	
	\begin{theorem}\label{thm4OminoLockedInfiniteConstruction}
		For any positive integer $n$, the square grid of side length $26 + 16n$ has a locked 4-omino tiling.
	\end{theorem}
	
	\begin{proof}
		The construction for $n = 1$, on the $42 \times 42$ grid, is presented in Figure \ref{fig4OminoLockedInfiniteConstructionSmall}, and the construction for $n = 4$ is presented in Figure \ref{fig4OminoLockedInfiniteConstructionBig}. It is self-evident that this construction can extended to arbitrary $n$ by adding additional repetitions of the patterns on the border and the 8 ``spokes'' leading into the center. Each repetition increases the grid dimensions by 8 on each side, for a total grid dimension increase of 16.
	\end{proof}
	
	\ipnc{.075}{T90-02.drawio}{\label{fig4OminoLockedInfiniteConstructionBig}A locked 4-omino tiling on the $90 \times 90$ grid. This is the tiling for $n = 4$ in the construction from Theorem \ref{thm4OminoLockedInfiniteConstruction}.}
	
	\subsection{A Locked 6-Omino Tiling}\label{sub6}
	
	A highly-optimized version of this search procedure, written in C++ by ChatGPT, was able to discover a locked 6-omino tilings on a $48 \times 48$ grid. We do not claim that this is the smallest example (as with the 4- and 5-omino tilings in Figures~\ref{figUniqueTilings10} and \ref{figUniqueTilings20}), nor are we certain that this is the only one on a $48 \times 48$ grid. The tiling is shown in Figure~\ref{figSampleTiling48}.
	
	\subsection{No Locked 7-Omino Tilings}\label{sub7}
	
	Our final computational result is the following impossibility:
	
	\begin{theorem}\label{thmNoLocked7OminoTilings}
		For any positive integers $m$ and $n$, the $m \times n$ grid does not have any locked 7-omino tilings.
	\end{theorem}
	
	The proof proceeds by a (gigantic, computationally-intensive) case analysis of what the tiling must look like near one of the corners. This result was also discovered by the ChatGPT-written algorithm, so it has been formalized in Lean. Below are some of the key definitions and the formal theorem statement; terms like \texttt{Cell}, \texttt{Connected}, and \texttt{tiles} have the obvious meanings. The full Lean project is available online.\footnote{\iftagged{name}{\url{https://github.com/jtuckerfoltz/LockedPolyominoTilings/tree/main/Lean}}{Link excluded for the purposes of anonymity.}}
	
	\begin{figure}
	\begin{lstlisting}
		def RigidPair (size : Nat) (a b : Finset Cell) : Prop :=
		∀ c d : Finset Cell, Connected c → Connected d → c.card = size → d.card = size →
		Disjoint c d → c ∪ d = a ∪ b → (c = a ∧ d = b) ∨ (c = b ∧ d = a)
		
		def Locked {size height width : Nat} (t : Tiling size height width) : Prop :=
		∀ a ∈ t.tiles, ∀ b ∈ t.tiles, a ≠ b → RigidPair size a b
		
		def Nondegenerate {size height width : Nat} (t : Tiling size height width) : Prop :=
		∃ u : Tiling size height width, u.tiles ≠ t.tiles
		
		theorem no_locked_nondegenerate_rectangle {height width : Nat}:
		¬ ∃ t : Tiling 7 height width, Nondegenerate t ∧ Locked t
	\end{lstlisting}
	\end{figure}	
	
	\ipnc{.16}{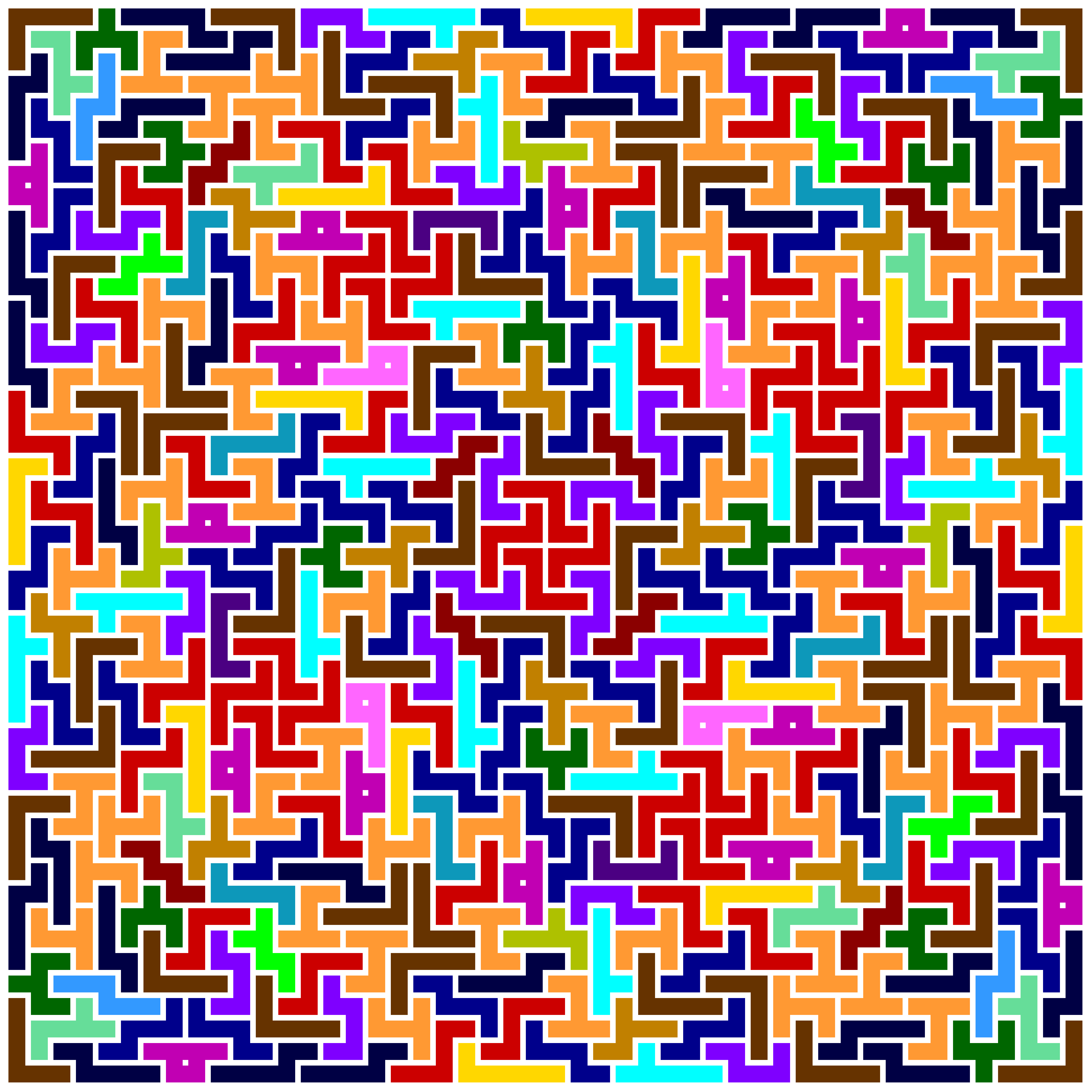}{\label{figSampleTiling48}The only locked 6-omino tiling known to exist.}
	
	\section{Infinite Tilings}\label{secTori}
	
	As evidenced by the proof of Theorem~\ref{thmNoLocked7OminoTilings}, a central difficulty in constructing locked polyomino tilings is making straight borders and corners. In this section, we remove this obstacle by considering infinite grid tilings. Any infinite tiling that is periodic can also be thought of as a tiling of a finite torus graph $T(m, n)$, i.e., the grid graph $G(m, n)$ with additional edges connecting the vertices on the left border of the grid to those on the right border and similarly with the top and bottom vertices. As Figure \ref{figTorusTilings} illustrates, with infinite grids we can construct several different kinds of tilings that are not possible for finite grids.
	
	\ipnc{.16}{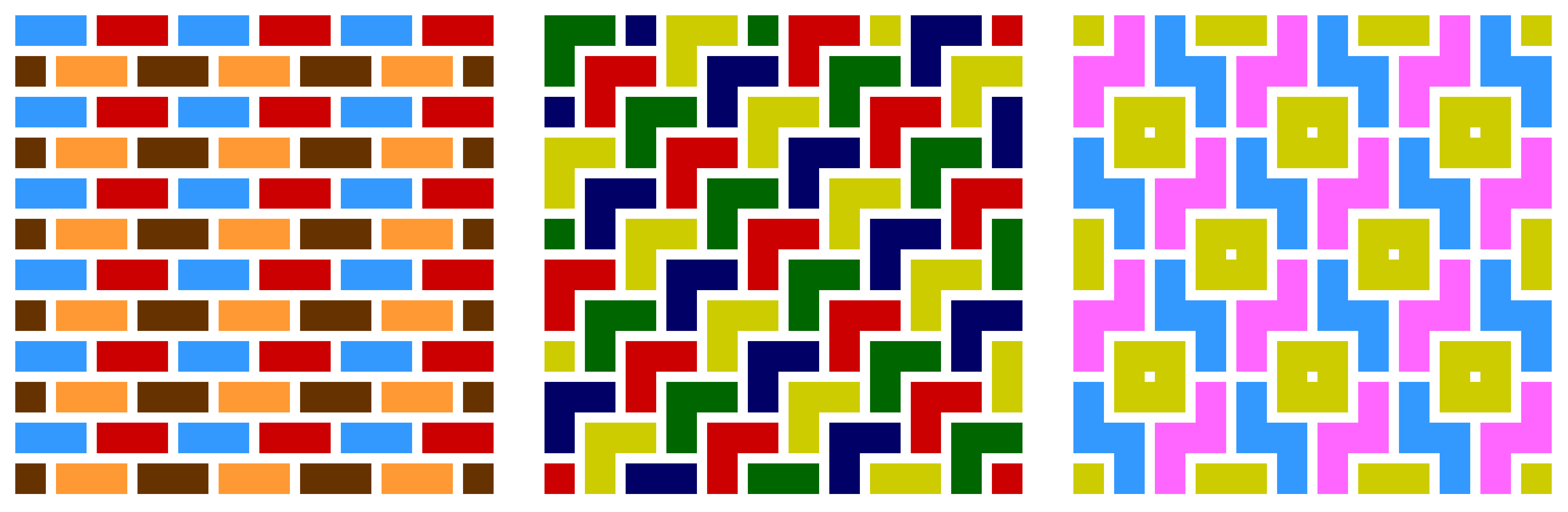}{\label{figTorusTilings} Some locked 2-, 3-, and 4-omino tilings on $T(12, 12)$ or the infinite grid. Note that the middle tiling is locked even if we were to allow recombining up to four tiles at a time rather than two. Generalizations of the tiling on the left are still locked with respect to arbitrarily large numbers of tiles being recombined.}
	
	Surprisingly, we find that it is possible to construct infinite periodic locked $t$-omino tilings for \emph{arbitrarily large} $t$. We can do this in either the square grid lattice, which is the focus of this paper, or the triangular lattice, which serves as an intriguing counterpoint to the connectivity result of Cannon \cite{TriangularLattice}. Both constructions are essentially the same, and we only sketch them here by means of example.
	
	\ipnc{.32}{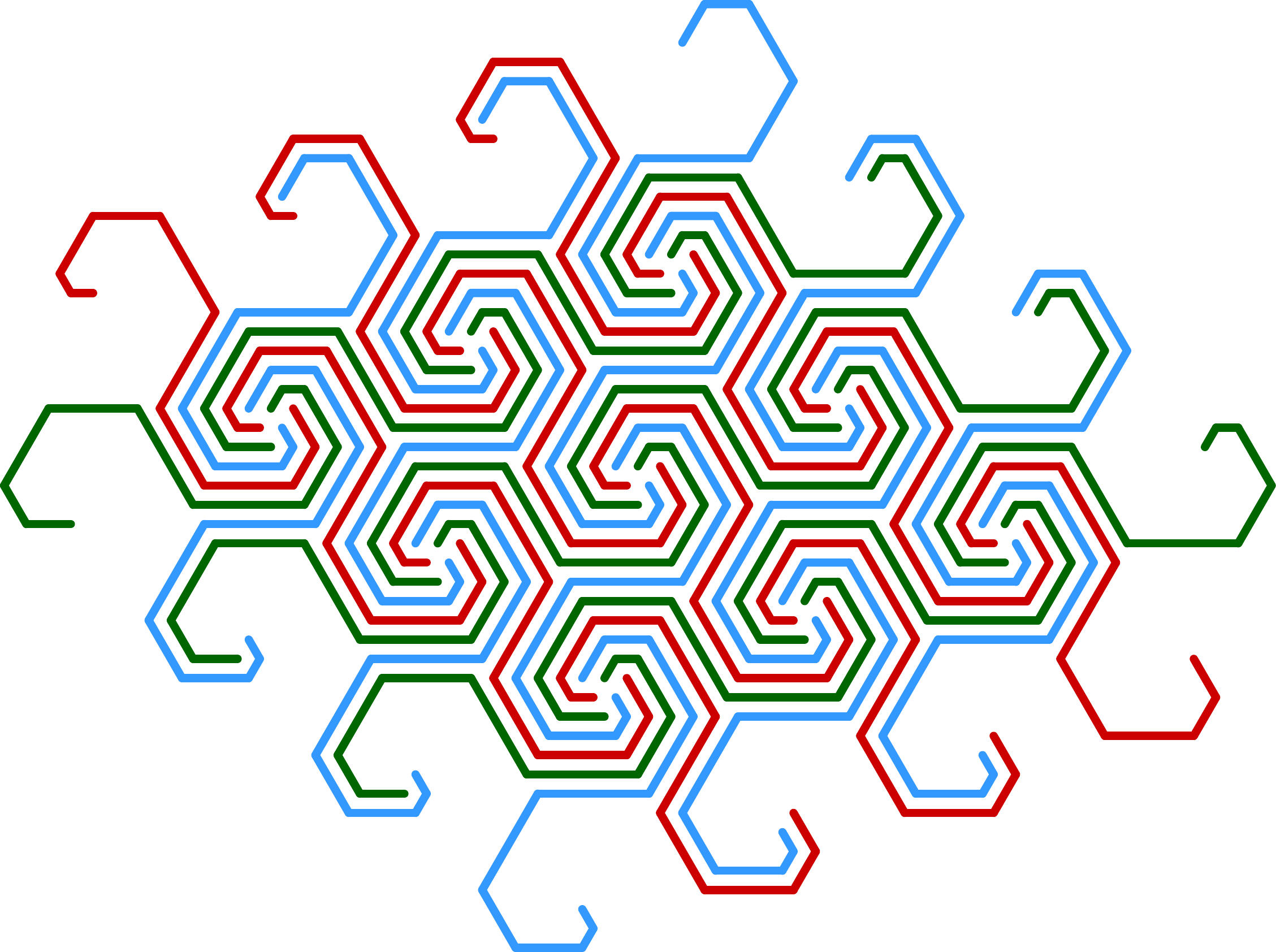}{\label{figBigTriangularTorusTiling} Part of a periodic locked 31-omino tiling of the infinite triangular lattice.}
	
	\begin{figure}[H]
		\centering
		\includegraphics[scale=.099]{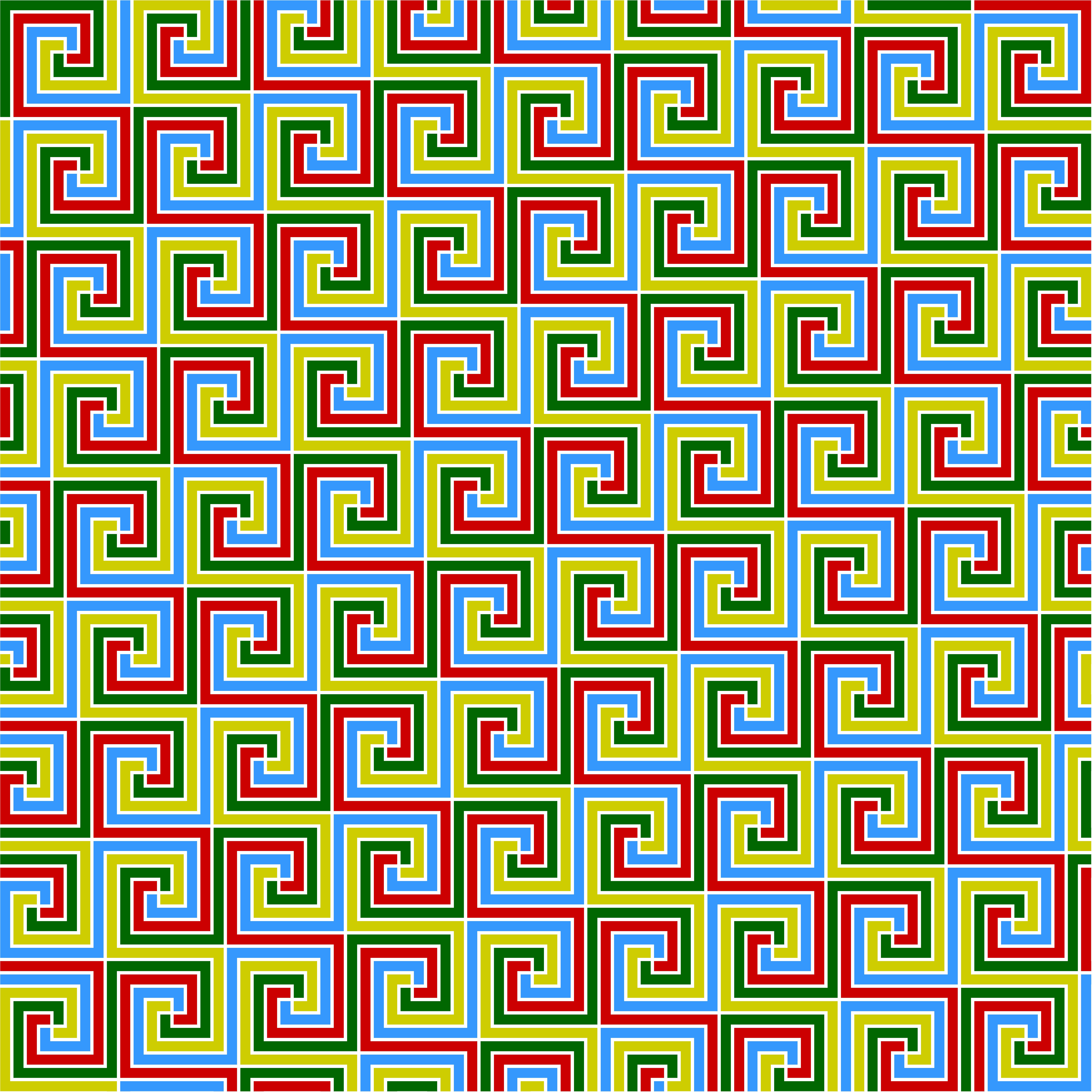}
		\caption{\label{figBigTorusTiling} A locked 41-omino tiling of the infinite grid, or the $82 \times 82$ torus.}
	\end{figure}
	
	\begin{theorem}\label{thmTorusFamily}
		For any positive integer $n$, the square torus $T(2(2n^2 + 2n + 1), 2(2n^2 + 2n + 1))$ (or infinite square grid lattice) has a locked $(2n^2 + 2n + 1)$-omino tiling, and the infinite triangular lattice has a locked $((n + 1)(n + 2) + 1)$-omino tiling.
	\end{theorem}
	
	\begin{proof}
		The construction is tedious to specify formally, and it is easy to verify that it generalizes to any $n$ from seeing a sufficiently large example, so we will mostly rely on illustrations. For the square lattice, the construction for $n = 1$ is the 5-omino tiling pattern in the center of Figure \ref{figSpaceInvaders}. The Z-tiles in this tiling can be thought of as double spirals, where we start with a straight segment of length 3 and spiral inward on both sides by appending paths of length 1 to the sides in a symmetric way. For general $n$, we extend this double spiral pattern in the natural way: we start with a straight segment of length $2n + 1$ and append symmetric spirals with segments of length $2n - 1, 2n - 3, 2n - 5, \dots, 3, 1$. Figure \ref{figBigTorusTiling} shows the tiling for $n = 4$. Every pair of adjacent tiles touches in the exact same way, so it suffices to verify that they cannot be recombined. This follows from the observation that the lengths of the two long paths trailing off of the region obtained by removing two adjacent tiles sum to $t + 1$, so they cannot be in the same $t$-omino. This forces them to partition up the rest of the region in a unique way. Figure \ref{figBigTriangularTorusTiling} shows the tiling for the triangular lattice (also for $n = 4$). It is straightforward to see that both patterns extend to arbitrary $n$ by making the spirals longer.
	\end{proof}
	
	Finally, we note that the metagraph of 4-omino infinite grid tilings does not merely consist of isolated vertices, but also larger connected components. This is intriguing because it is the only known result of its kind for $t = 4$.
	
	\begin{proposition}
		For any positive integers $m$ and $n$, $\mathcal{M}(T(4m, 5 + 8n), 4)$ has connected components of size $3^m$. Consequently, the metagraph of the infinite grid, $\mathcal{M}(\zz \times \zz, 4)$, has multiple connected components of infinite size.
	\end{proposition}
	
	\ipnc{.16}{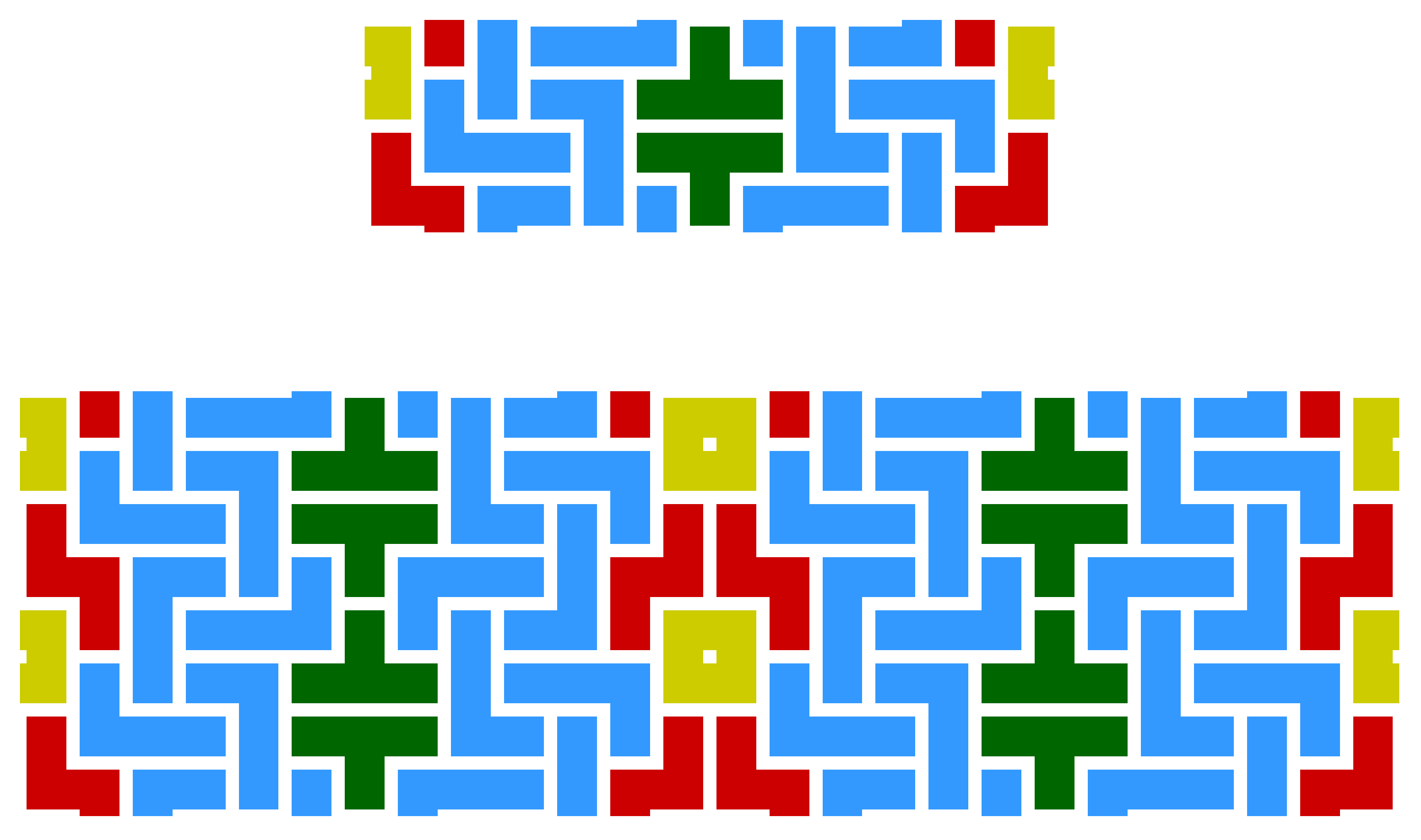}{\label{figTorusComponent}Above, a tiling in a connected component of size 3 in $\mathcal{M}(T(4, 13), 4)$ that can be tiled with itself to produce arbitrarily large connected components in torus metagraphs. Below, four copies are joined together, illustrating how it would extend to an infinite grid.}
	
	\begin{proof}
		Consider the torus tiling in Figure \ref{figTorusComponent}. This can be extended to a tiling of $T(4m, 5 + 8n)$ by placing copies of the blue L-tiles a total of $n$ times on each side, then tiling the whole diagram with itself vertically for a total of $m$ copies. Only the green T-tiles are unlocked. Each pair of T-tiles may be independently re-tiled in 3 different ways (the original T-tiling plus two different Z-tilings), but none of these admit recombination moves with neighboring tiles.
	\end{proof}

	\section{Weighted Grid Tilings}\label{secSixDistricts}
	
	In this section we consider an extension where vertices may have different weights, as in the motivating setting of political redistricting, where the weights are populations. Specifically, we consider the minimal extension where each vertex has weight 1 or 2. To avoid having to introduce new notation for this setting, we can model weight-2 vertices by appending leaves to each of them. Formally, we define an $m \times n$ \emph{weighted grid graph} to be a graph obtained by starting from an $m \times n$ grid graph $G(m, n)$ and adding some number of new degree-one vertices, each adjoined to distinct vertices from $G(m, n)$. Our main result is that weighted grid graphs admit arbitrarily large locked tilings with a constant number of tiles; an example of such an instance is depicted in Figure~\ref{figSixDistrictsInstance201}.
	
	\ignc{.37}{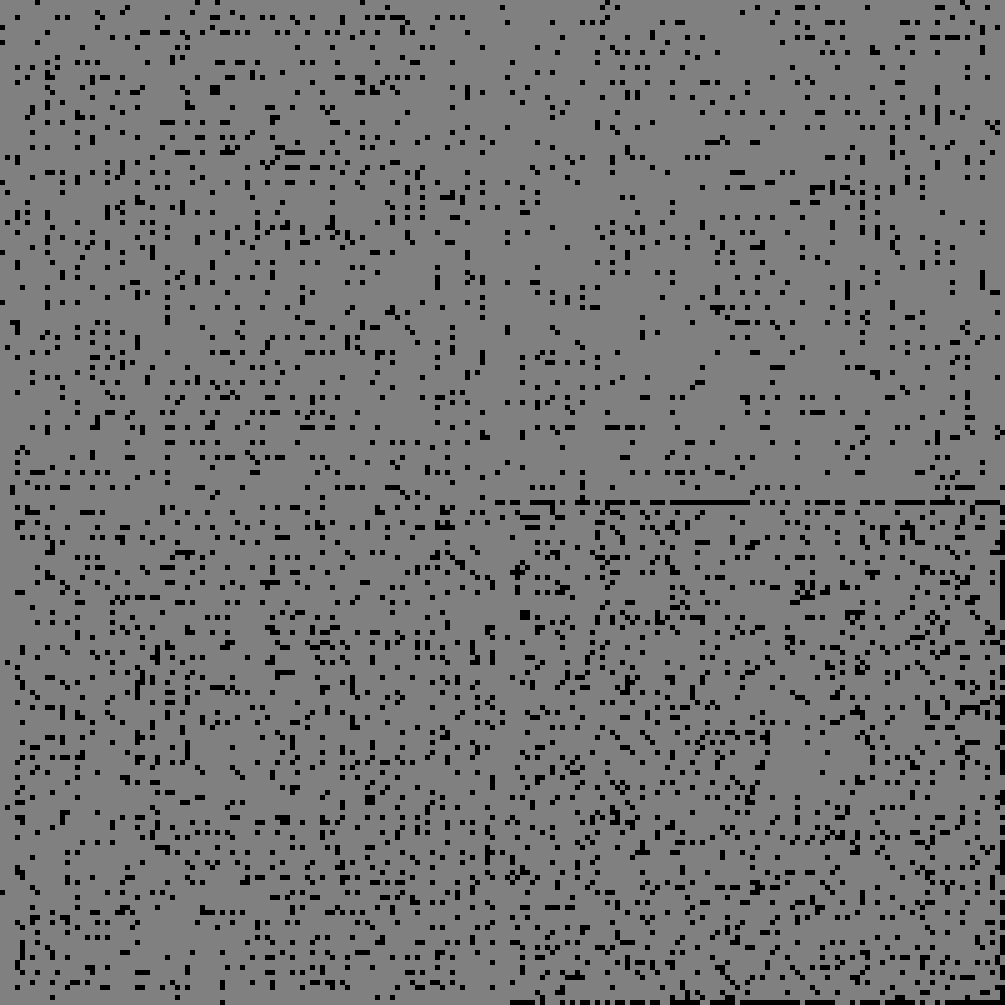}{\label{figSixDistrictsInstance201}A $201 \times 201$ weighted grid graph on which the 6-district metagraph is disconnected. Grey pixels represent vertices of weight 1, and black pixels represent vertices of weight 2.}
	
	\begin{theorem}\label{thmSixDistricts}
		For any odd integer $n \geq 53$, there exists an $n \times n$ weighted grid graph $G = (V, E)$ such that $\mathcal{M}(G, \frac{\abs{V}}{6})$ has an isolated vertex.
	\end{theorem}

	This result is striking because results about balanced connected partitions of grid graphs can often be extended to \emph{grid-like} graphs \cite{procacciaTuckerFoltz2021compactness, TreeSplitting}. Weighted grid graphs are special cases of both these definitions; they are aesthetically about as similar to grids as one can get. This result shows that Conjecture~\ref{cnjGridArbitraryK}, if true, does not extend in an analogous way. Even in the regime of small $k$ and large $n$, metagraph connectivity is an extremely brittle property that seems to require the uniformity and symmetry properties of grids.
	
\begin{figure}[ht]
		\centering
		\begin{subfigure}{.5\textwidth}
			\centering
			\includegraphics[scale=.33]{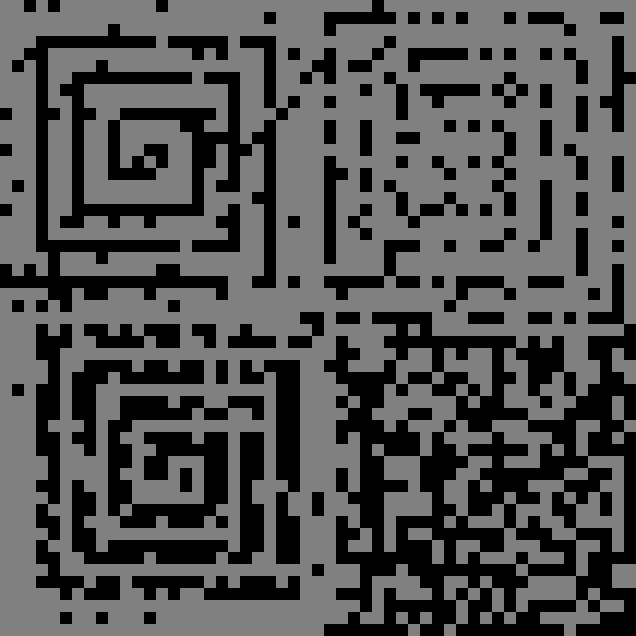}
		\end{subfigure}%
		\begin{subfigure}{.5\textwidth}
			\centering
			\includegraphics[scale=.1414]{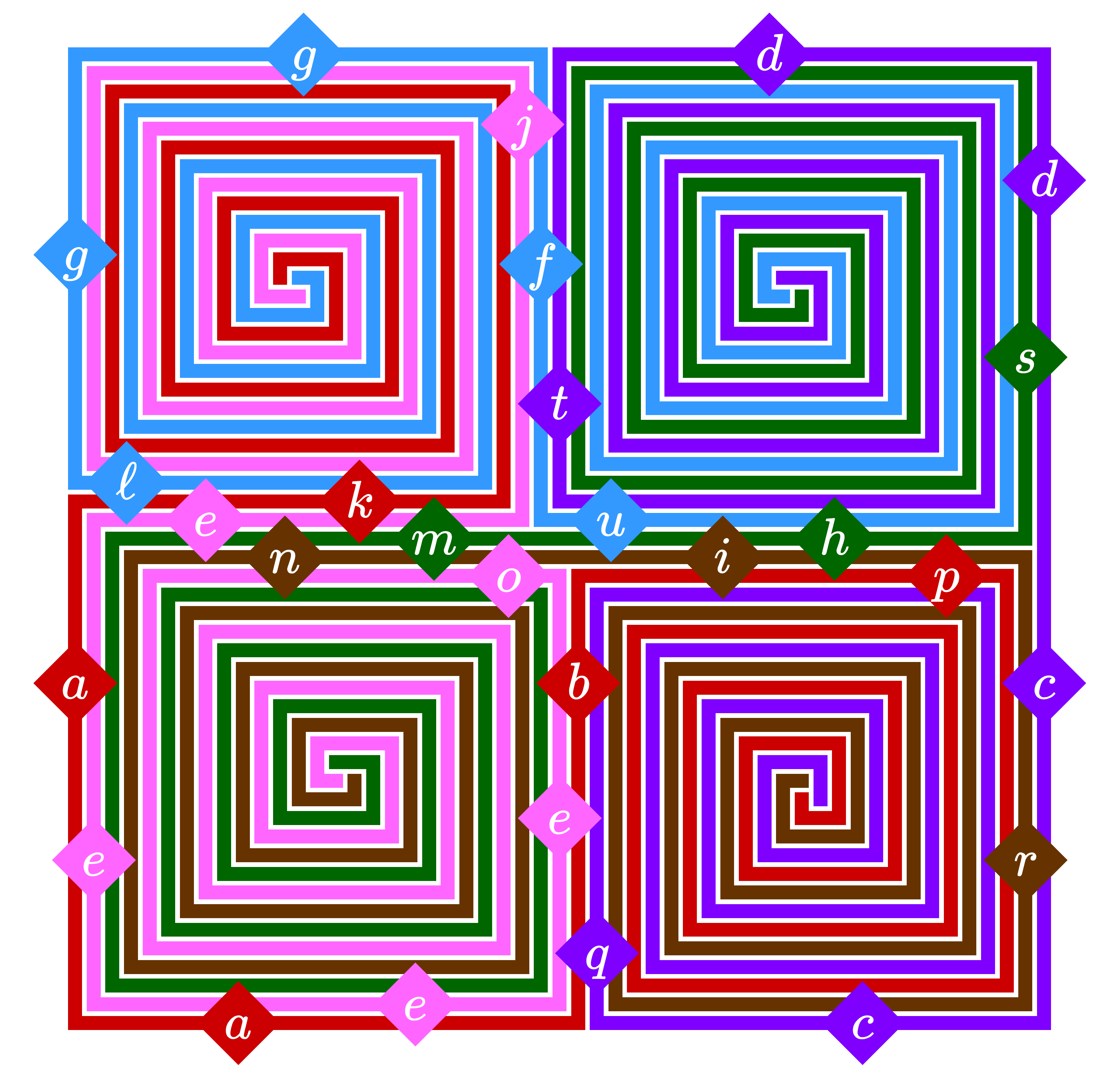}
		\end{subfigure}
		\caption{A $53 \times 53$ weighted grid graph (weight 2 vertices in black again) and a locked tiling on it. There are 1043 vertices of weight 2, for a total weight of 3852, divided evenly among 6 tiles, each of weight 642.}
		\label{figSixDistrictsInstance53}
	\end{figure}

	The proof follows from a general construction of locked tilings similar to Figure~\ref{figBigTorusTiling}. Every tile starts and ends in a spiral, and we take care to ensure that the entire diagram can be made square. Figure~\ref{figSixDistrictsInstance53} shows the smallest tiling from our construction, on a $53 \times 53$ grid; this structure generalizes to any odd grid side length $n$ by adding some number of rows and columns to each spiral. We will not be completely formal in arguing that the construction does indeed work for all odd $n \geq 53$, but we will give a description of how to find appropriate weightings of any sufficiently large grid.
	
	As in the constructions from Theorem~\ref{thmTorusFamily}, the tiling has the property that the union of every pair of adjacent tiles forms a very specific kind of shape: a corridor of width 2, terminating in a spiral, with two long paths trailing off the other end. So long as the sum of weights of every pair of trailing paths is greater than the target tile size, we have a locked tiling for the same reason as before. For sufficiently large $n$, this reduces to finding \emph{multipliers} for the populations of each of the segments labeled by diamonds in Figure~\ref{figSixDistrictsInstance53}. For example, a multiplier of $g = 1.5$ would mean that the total population of the $g$ segments in the top-left boundary should have $1.5 \times$ their ordinary population in the unweighted case, up to an additive error that does not grow with the grid dimensions. One could implement this by assigning a random half of the vertices in these segments to have weight 2 rather than 1. As long as the multipliers of all \emph{interior} segments (meaning not the ones going into a spiral) are all between 1 and 2, we will be able to find a valid assignment of weights as long as we let $n$ grow large enough. To find such a set of multipliers, we solve the following linear program, where variable $T$ represents the sum of multipliers across an entire tile.
	$$\lpmin{\max\{a, b, c, d, e, f, g, h, i\}}{\\
		
		$a, b, c, d, e, f, g, h, i \geq 1$ & $T \leq u + f + 2g     +   2a + b + p$\\
		$j, k, \ell, m, n, o, p, q, r, s, t, u \geq 0$ & $T \leq \ell + 2g         +   2d + 2c + q$\\
		& $T \leq t + 2d + 2c   +   2a + k$\\
		
		$T = \ell + 2g + f + u$ & $T \leq j               +   h + s$\\
		$T = j + 4e + o$ & $T \leq \ell + 2g + f     +   m$\\
		$T = k + 2a + b + p$ & $T \leq o + 4e         +   u$\\
		$T = s + h + m$ & $T \leq t + 2d         +   i + n$\\
		$T = t + 2d + 2c + q$ & $T \leq s               +   r$\\
		$T = n + i + r$ & $T \leq m + h           +   2c + q$\\
		& $T \leq o               +   p$\\
		& $T \leq n               +   b + 2a + k$\\
		& $T \leq j + 4e         +   i + r$\\
	}$$
	Here the equality constraints enforce that each tile has the same size, and the right column of inequality constraints enforce that every adjacent pair of tiles cannot be recombined, because the trailing paths are too large. We allow multipliers $j, k, \dots, u$ to be less than one, and do not care how large they get, because we will end up adding a large population to each of these pieces in each spiral anyway. The optimal solution sets $c = h = \frac53$ and $a = b = d = e = f = g = i = 1$. The key property of this optimal solution is that all values are strictly less than 2; had they been higher, we would have accordingly needed higher popoluation weights on vertices within the interior segments. It is clear in Figure~\ref{figSixDistrictsInstance201} that the $c$ and $h$ segments have indeed been up-weighted the most: $\frac23$ of the segments have weight 2, for an average weight of $\frac53$.
	
	This analysis only works in the limiting case as the grid side length $n$ is sufficiently large. For each finite $n$, the multipliers will need to be different by progressively tinier amounts. The first case where the actual multipliers are valid is $n = 53$. From there, the constraints only become easier to satisfy.
	
	\section{Open Problems}\label{secOpenProblems}
	
	The practical success of the ReCom sampling algorithm has spurred renewed interest in the connectivity of polyomino tilings under local moves. We are still quite far from being able to prove Conjectures \ref{cnjGridKEquals3} and \ref{cnjGridArbitraryK}, and this paper hopes to have shed light on some of the obstructions that any potential proof will need to overcome. We have demonstrated that, even in the tame setting of grid graphs, certain local structures of balanced graph partitions exist which can disconnect the state space of ReCom.
	
	The study of these local structures is quite intriguing on its own, due to the fact that locked polyomino tilings are so difficult to find. Numerous open questions remain, of which the following are particularly interesting.
	
	\begin{enumerate}
		\item Are there locked 5- and 6-omino tilings for arbitrarily large grids? Can we construct an infinite family? None of the 3271 tilings enumerated in this paper appear to generalize.
		\item Are there any locked $t$-omino grid tilings, of any finite grid dimensions, for any $t > 7$? Is the set of positive integers $t$ with this property bounded?
		\item Can locked $t$-omino tilings exist even if we allow 3-way recombination moves? For infinite grids/tori, the answer is clearly yes, as shown in Figure \ref{figTorusTilings}. For finite grids, however, this is open for $t \geq 4$. (For $t = 3$, a simple enumeration of all cases for what could be placed in a corner shows that there are no locked tilings.) If we allow 4-way recombination moves, it is unknown whether locked tilings exist even on tori.
		\item What structures exist among locked 3-omino tilings, and more generally, 3-omino grid metagraphs? As it so happens, 3-omino metagraphs suffer from more connectivity issues than just isolated points. The connected components of $\mathcal{M}(G(6, 6), 3)$ are enumerated in Table \ref{tabComponents}.
		
		\begin{table}[H]\begin{center}
				\begin{tabular}{ c|c }
					Number of vertices in component & Number of components\\
					\hline
					73738 & 1\\
					384 & 4\\
					235 & 8\\
					199 & 8\\
					68 & 8\\
					20 & 16\\
					19 & 16\\
					16 & 2\\
					8 & 8\\
					2 & 16\\
					1 & 50\\
				\end{tabular}
				\caption[]{\label{tabComponents}Connected component sizes of $\mathcal{M}(G(6, 6), 3)$.}
		\end{center}\end{table}
		It is not surprising that there are as many as 50 locked tilings (components of size 1) or that most of the vertices lie in one giant component. What \emph{is} surprising is that there are many other large components. For all components of size 68 and above, there is no pair of grid cells that are always together in the same 3-omino throughout the component. In other words, some components do not correspond to local regions of the $6 \times 6$ grid, but rather global invariants that are preserved through recombination moves. Figure~\ref{figGraphVis} depicts one of these strange components.
		
		\ignc{.29}{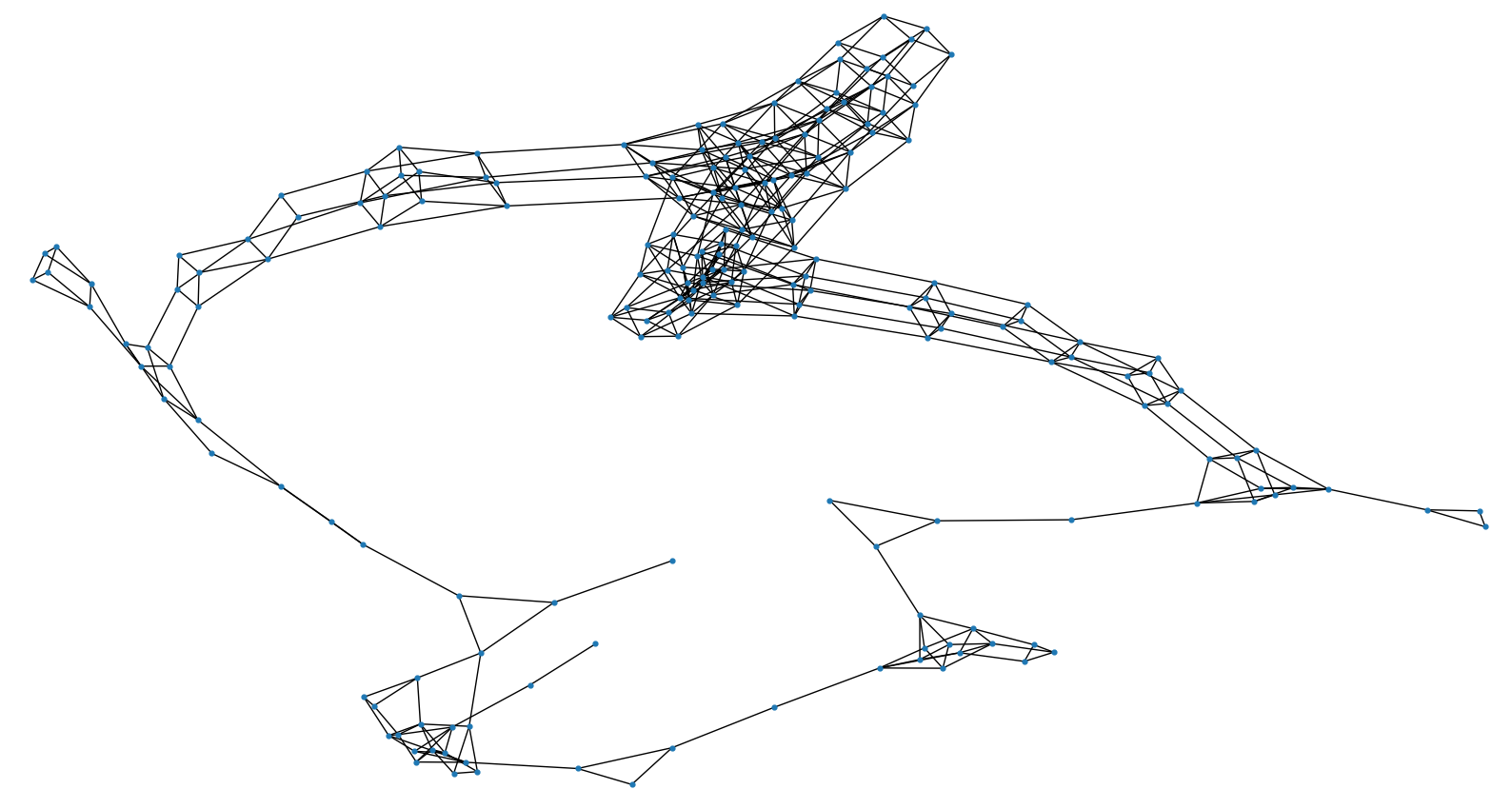}{\label{figGraphVis} A connected component of $\mathcal{M}(G(6, 6), 3)$ of order 199. Every vertex represents a 3-omino tiling of the $6 \times 6$ grid. The embedding in the plane is arbitrary, but it illustrates that the graph has the ``topology'' of a circle in an informal sense.}
		
		A forthcoming work by Andrew Chen, Sara Stephens, Gabe Udell, and Jonathan Webb provides an invariant that partially explains this phenomenon: These components all have different invariant values from the main component. However, there are still multiple large components that have the same invariant values as each other, so there may still be some hidden structure to discover.
		\item Even when a grid metagraph is disconnected, is it still the case that the vast majority of tilings lie in one connected component? The result mentioned in the previous paragraph indicates that the answer is likely negative for $n = 3$. However, the rest of the results in this paper indicate the answer is likely affirmative for $n \geq 4$.
	\end{enumerate}

	\tagged{name}{
		\section*{Acknowledgments}
		
		I am very grateful to Christopher Donnay for his helpful suggestions, motivations, pointers to existing literature, and technical support. I am also grateful to Parker Rule, Moon Duchin, and Josh Mermelstein for multiple thought-provoking discussions. And also to Peep Mouse Kumar, the kind and fluffy cat that sat on my lap (and sometimes arms) as I wrote the first draft of this paper.
		
		This material is based in part upon work supported by the National Science Foundation under Grant No.\txt{} DMS-1928930 and by the Alfred P.\txt{} Sloan Foundation under grant G-2021-16778, while the author was in residence at the Simons Laufer Mathematical Sciences Institute (formerly MSRI) in Berkeley, California, during the Fall 2023 semester.}
		
		The results of Sections \ref{sub6} and \ref{sub7} were obtained using GPT-6 Astra.
	\newpage
	
	\bibliographystyle{plain}
	\bibliography{bibliography}
\end{document}